\documentstyle[amscd,amssymb,verbatim,diagrams,12pt]{amsart}
\newcommand{\hH}{{\Bbb H}}

\newcommand{\Sets}{\operatorname{Sets}}

\newcommand{\ZZ}{{\cal Z}}
\newcommand{\Id}{\operatorname{Id}}

\renewcommand{\mod}{\operatorname{mod}}
\newcommand{\Com}{\operatorname{Com}}

\newcommand{\gr}{\operatorname{gr}}

\newcommand{\und}{\underline}

\newcommand{\OO}{{\cal O}}

\newcommand{\coker}{\operatorname{coker}}

\newcommand{\NN}{{\cal N}}

\newcommand{\II}{{\cal I}}

\newcommand{\hra}{\hookrightarrow}

\newcommand{\GG}{{\cal G}}
\newcommand{\CC}{{\cal C}}
\newcommand{\UU}{{\cal U}}

\newcommand{\Spec}{\operatorname{Spec}}

\newcommand{\Mat}{\operatorname{Mat}}

\newcommand{\si}{\sigma}

\newcommand{\de}{\delta}

\newcommand{\im}{\operatorname{im}}

\newcommand{\A}{{\Bbb A}}

\numberwithin{equation}{subsection}

\newcommand{\GL}{\operatorname{GL}}

\newtheorem{thm}{Theorem}[subsection]
\newtheorem{prop}[thm]{Proposition}
\newtheorem{lem}[thm]{Lemma}

{  \theoremstyle{definition}
\newtheorem{defi}[thm]{Definition}
\newtheorem{ex}[thm]{Example}

\newtheorem{rem}[thm]{Remark}

}

\newcommand{\Pf}{\noindent {\it Proof}}
\newcommand{\id}{\operatorname{id}}
\newcommand{\Lie}{\operatorname{Lie}}
\newcommand{\ov}{\overline}
\newcommand{\we}{\wedge}

\newcommand{\Ad}{\operatorname{Ad}}

\renewcommand{\AA}{{\cal A}}

\newcommand{\FF}{{\cal F}}
\newcommand{\EE}{{\cal E}}

\newcommand{\TT}{{\cal T}}

\newcommand{\HH}{{\cal H}}

\newcommand{\VV}{{\cal V}}
\newcommand{\SS}{{\cal S}}

\newcommand{\Om}{\Omega}

\newcommand{\Hom}{\operatorname{Hom}}

\newcommand{\Ext}{\operatorname{Ext}}
\newcommand{\End}{\operatorname{End}}

\newcommand{\Aut}{\operatorname{Aut}}

\newcommand{\Eend}{{\cal E}nd}
\newcommand{\Hhom}{{\cal H}om}
\renewcommand{\a}{\alpha}

\newcommand{\om}{\omega}
\newcommand{\De}{\Delta}

\newcommand{\th}{\theta}
\newcommand{\C}{{\Bbb C}}

\newcommand{\La}{\Lambda}
\newcommand{\Ga}{\Gamma}

\newcommand{\wt}{\widetilde}
\newcommand{\ot}{\otimes}

\newcommand{\sub}{\subset}
\newcommand{\ed}{\qed\vspace{3mm}}

\newcommand{\Fun}{\operatorname{Fun}}

\newcommand{\car}{\curvearrowright}

\newcommand{\Der}{\operatorname{Der}}

\newcommand{\I}{{\mathcal{I}}}
\newcommand{\LLie}{{\cal L}ie}
\newcommand{\Isom}{\operatorname{Isom}}
\newcommand{\Inn}{\operatorname{InnAut}}

\title{NC-smooth algebroid thickenings for families of vector bundles and quiver representations}
\author{Ben Dyer}
\address{BD: The Evergreen State College, Olympia, WA 98505} 
\author{Alexander Polishchuk}
\address{AP: University of Oregon and National Research University Higher School of Economics}

\begin{document}

\begin{abstract}
In his work on deformation quantization of algebraic varieties Kontsevich introduced the notion of {\it algebroid}
as a certain generalization of a sheaf of algebras. 
We construct algebroids which are given locally by NC-smooth thickenings
in the sense of Kapranov, over two classes of smooth varieties:
the bases of miniversal families of vector bundles on projective curves,
and the bases of miniversal families of quiver representations.
\end{abstract}

\maketitle

\section*{Introducton}

In this work we study certain families of vector bundles over noncommutative bases.
More precisely, our framework is the theory of NC-schemes over $\C$, developed by Kapranov in \cite{Kapranov}.
These are analogs of usual schemes based on the algebras that are close to being commutative: any expression containing
sufficiently many commutators in such rings vanish. More precisely, these are {\it NC-nilpotent} algebras; one also considers {\it NC-complete} algebras which are complete with respect to the commutator filtration.

In this theory there is a natural notion of {\it NC-smoothness}, which is analogous to the notion of quasi-free algebra
from \cite{CQ}. Kapranov proves the existence and uniqueness
of an NC-smooth thickening for any smooth affine scheme $X$. By definition, such a thickening corresponds to an NC-smooth
algebra whose abelianization is the algebra of functions on $X$.
The problem of determining which non-affine
smooth schemes admit such extensions seems to be quite hard. There are very few known examples of such thickenings.
For example, there are explicit constructions for Grassmannians and abelian varieties (see \cite{Kapranov},\cite{PT}). In both cases the relevant NC-smooth thickenings
represent natural functors on the category $\NN$ of NC-nilpotent algebras.
On the other hand, there is no smooth scheme for which we would know that there is no NC-smooth thickening.

One of the constructions considered in \cite{Kapranov} is that of a natural functor 
of families of vector bundles over NC-nilpotent bases, which on the commutative level are induced by a given 
family of vector bundles on a fixed projective variety with a base $B$. 
More precisely, we consider the following situation.
Let $Z$ be a projective algebraic variety, $B$ a smooth variety, and let $\EE^{ab}$ be a vector bundle over $B$.
We denote by $\rho: B\times Z \to B$ the natural projection. 

\begin{defi}\label{versal-fam-defi}
We say that $\EE^{ab}$ is an {\it excellent family} of bundles on $Z$ if
\begin{enumerate}
	\item[(a)] $\OO_B \to \rho_*\Eend(\EE^{ab})$ is an isomorphism,
	\item[(b)] the Kodaira-Spencer map $\kappa: T_B \to R^1\rho_*\Eend(\EE^{ab})$ is an isomorphism,
	\item[(c)] $R^2\rho_*\Eend(\EE^{ab}) = 0$,
	\item[(d)] $R^d\rho_*\Eend(\EE^{ab})$ is locally free for $d \geq 3$.
\end{enumerate}
\end{defi}

For example, if $Z$ is a projective curve then conditions (c) and (d) are automatic. Condition (a) is satisfied
for a family of simple bundles (see \cite[Lem.\  4.6.3]{HL}). Condition (b) is satisfied if the map from $B$ to
the moduli stack of vector bundles on $Z$ is \'etale.

Following \cite{Kapranov} we consider the natural functor $h^{NC}_B$ on the category $\NN$ of noncommutative families of vector bundles compatible with $\EE^{ab}$ (see Definition \ref{fun-defi} for details).
It was claimed in \cite{Kapranov} that this functor is representable by an NC-smooth thickening of $B$.
However, the proof contained a gap. 

In the present paper we prove that whenever $\dim B\ge 1$, the functor $h^{NC}_B$ is not representable by an NC-scheme
(see Theorem \ref{non-rep-thm}). The reason for this is rather silly: we observe that $h^{NC}_B$ factors through
the quotient category $a\NN$ of $\NN$ in which conjugate homomorphisms are identified (see Sec.\ \ref{aNC-sec}). 

The natural idea then is to ask the representability question in this new category $a\NN$.
Our main technical result is that this is true locally: 
the functor of families over NC-nilpotent bases is representable in the case when $B$ is affine
(see Theorem \ref{loc-rep-thm}). 
We use this local representability of $h^{NC}_B$ in $a\NN$ to construct in the general case a 
{\it $\C$-algebroid}\footnote{This notion has nothing to do with the more commonly used {\it Lie algebroid}: the latter
is a sheaf of Lie algebras with some extra structures, whereas a $\C$-algebroid is a certain stack of $\C$-linear
categories.}
 over $B$ in the sense of \cite{Kont}, \cite[Sec.\ 2.1]{KS}, given locally
by an NC-smooth thickening of $B$. We call such a structure an {\it NC-smooth algebroid thickening} of $B$ 
(see Definition \ref{NC-stack-alg-defi} for details). 

\medskip

\noindent
{\bf Theorem A}[see Thm.\ \ref{functor-algebroid-thm}+Thm.\ \ref{loc-rep-thm}]. {\it Let $B$ be a (smooth) base of an 
excellent family of vector bundles.
Then there exists an NC-smooth algebroid thickening of $B$. 
}

\medskip

In the case when $B$ is quasi-projective, so that there exists an
open affine covering $(U_i)$ of $B$, such that all intersections $U_i\cap U_j$ are distinguished affine opens in both $U_i$
and $U_j$, the algebroid in Theorem A can be described in more down-to-earth terms as follows.
We have an NC-smooth thickening of $U_i$ for each $i$; over $U_i\cap U_j$ we have isomorphisms between 
the two induced NC-smooth thickenings; and over $U_i\cap U_j\cap U_k$ the isomorphisms agree up to an inner automorphism
(furthermore, the corresponding invertible elements are chosen and satisfy the natural compatibility condition over 
$U_i\cap U_j\cap U_k\cap U_l$).

Note that algebroids were introduced by Kontsevich in connection with deformation quantization of algebraic varieties (see \cite{Kont}, \cite{KS}). NC-smooth thickenings are in some ways similar to deformation quantization algebras (in particular,
the construction of NC-smooth thickenings from torsion-free connections in \cite{PT} is somewhat reminescent of Fedosov's deformation quantization procedure in \cite{Fedosov}).
Thus, it is not very surprising that algebroids made their appearance in the theory of NC-smooth thickenings.
In light of Theorem A, it seems that rather than asking which smooth schemes admit NC-smooth thickenings,
it is more natural to ask which smooth schemes admit  NC-smooth algebroid thickenings.

In fact, in the proof of Theorem A we construct a canonical algebroid (up to an equivalence). 
Since there is a well defined notion of a module over an algebroid, one natural problem is whether there is a universal
family of modules over our algebroid, extending the original family over $B$. 
One can also try to study the higher rank analog of the Fourier-Mukai transform picture for NC Jacobians considered in
\cite[Sec.\ 4]{PT}. We leave these questions for a future study.

Motivated by Toda's work \cite{Toda}, we also consider the similar picture for representations of quivers. 
Namely, starting with an {\it excellent family} (see Def.\ \ref{versal-quiver-def}) of
representations of a finite quiver $Q$ (with no relations), we consider the functor of compatible families of representations
of $Q$ over NC-nilpotent affine schemes. We show that the situation is completely similar (and somewhat easier) to the 
picture discussed above. 

\medskip

\noindent
{\bf Theorem B}[see Thm.\ \ref{functor-algebroid-thm}+Thm.\ \ref{loc-rep-quiver-thm}]. {\it Let $B$ be a (smooth) base of an 
excellent family of representations of $Q$.
Then there exists an NC-smooth algebroid thickening of $B$. 
}

\medskip

For example, this result applies to the moduli space of stable quiver representations corresponding to an indivisible
dimension vector.

Note that for the proof of Theorem B we develop a version of nonabelian hypercohomology $\hH^1$ for a sheaf of groups
acting on a sheaf of sets, which may be of independent interest (see Section \ref{nonab-hypercoh}).

Toda also constructs in \cite{Toda} local (non NC smooth) NC thickenings for some obstructed families of
vector bundles (and for representations of quivers with relations). It would be interesting to study whether
these thickenings glue into an algebroid.



The paper is organized as follows.
In Section \ref{aNC-nonrep-sec} we discuss the category $a\NN$
of affine almost NC schemes (in which conjugate homorphisms are identified).
We prove in Section \ref{gluing-sec} that any formally smooth functor on $a\NN$, that is locally representable, leads to an
NC smooth algebroid thickening (see Theorem \ref{functor-algebroid-thm}). Then in Section \ref{versal-sec}
we show that the functor of NC families extending the given excellent family of vector bundles factors through $a\NN$,
and as a consequence, is not representable except in trivial cases (see Theorem \ref{non-rep-thm}).

In Section \ref{rep-sec} we prove local representability results for formally smooth 
functors on $a\NN$. 
First, we give a technical representability criterion
for such a functor extending the functor on commutative algebras representable by a smooth affine scheme
(see Proposition \ref{aNC-rep-prop}). Then we apply this criterion to the functor of NC families extending a given
excellent family of vector bundles (see Theorem \ref{loc-rep-thm}) and then to the functor of NC families  of 
quiver representations (see Theorem \ref{loc-rep-quiver-thm}).

\medskip

\noindent
{\it Acknowledgments}. The work of the second author is supported in part by the NSF grant DMS-1700642 and by the 
Russian Academic Excellence Project `5-100'.
He also would like to thank Institut de Math\'ematique de Jussieu and
Institut des Hautes \'Etudes Scientifiques for hospitality and excellent working conditions.

\medskip

\noindent
{\it Conventions}. All algebras we consider are over $\C$, all schemes are assumed to be of finite type over $\C$.
The expression $[a,b]$ always denotes commutator in an associative algebra: $[a,b]=ab-ba$.

\section{Affine almost NC schemes and the non-representability of the functor of NC families of vector bundles}
\label{aNC-nonrep-sec}

\subsection{Generalities on NC schemes}

For a ring $R$, we define the decreasing filtration $\I_nR$ by
$$
\I_nR=\sum_{i_1\ge 2,\ldots,i_m\ge 2,i_1+\ldots+i_m\ge n} R\cdot R_{i_1}^{\Lie}\cdot R\cdot\ldots\cdot R\cdot R^{\Lie}_{i_m}\cdot R,
$$
where $R_n^{\Lie}$ is the $n$th term of the lower central series of $R$ viewed as a Lie algebra.
Note that $R/\I_2R$ is precisely $R^{ab}$, the abelianization of $R$.

We define the category $\NN_d$ of {\it NC-nilpotent algebras of degree $d$} as the category of algebras $R$ 
for which $\I_{d+2}R=0$.
Thus, $\NN_0=\Com$ is the category of commutative algebras. A ring $R$ is in $\NN_1$ if and only if it is a
central extension of a commutative algebra. Here we say that an extension of algebras
$$0\to I\to R\to \ov{R}\to 0$$
is a {\it central extension} if $I$ is a central ideal in $R$ with $I^2=0$.

We denote by $\NN=\cup_{d\ge 0}\NN_d$ the category of NC-nilpotent algebras. 
For $A\in\NN$ we denote by $h_A$ the corresponding
representable covariant functor on $\NN$: $h_A(B)=\Hom_{alg}(A,B)$.

An algebra $R$ is called {\it NC-complete} if it is complete with respect to the filtration $(\I_nR)$.
We denote by $\NN\CC$ the category of NC-complete algebras.
For an NC-complete algebra $R$ we denote by $h_R$ the functor on $\NN$ given by $h_R(B)=\Hom_{alg}(R,B)$.
Note that the restriction $h_R|_{\NN_d}$ is naturally isomorphic to the representable functor $h_{R/\I_{d+2}R}$.
This easily implies that the functor 
$$\NN\CC^{op}\to \Fun(\NN,\Sets): R\mapsto h_R$$
is fully faithful.

An NC-complete algebra $R$ is called {\it NC-smooth} if the functor $h_R$ is formally smooth,
i.e., for any central extension in $\NN$, $B'\to B$, the induced map $h_R(B')\to h_R(B)$ is surjective.
An NC-nilpotent algebra $A$ of degree $d$ is called {\it $d$-smooth} if the same is true for the functor
$h_A|_{\NN_d}$.

Kapranov defines NC-nilpotent schemes (over $\C$) 
as locally ringed spaces locally isomorphic to the spectrum
of an NC-nilpotent algebra, with its natural structure sheaf, which is defined similarly to the commutative case.
General NC-schemes are similarly modeled on formal spectra of NC-complete algebras
(see \cite[Sec.\ 2]{Kapranov} for details). One can view an NC-scheme $X$ as an underlying usual scheme
$X^{ab}$ equipped with a sheaf of noncommutative algebras $\OO_X$ such that its abelianization is
$\OO_{X^{ab}}$. In this case we say that $X$ is an {\it NC-thickening of} $X^{ab}$. In the case when $X$
is NC-smooth, we say that it is an {\it NC-smooth thickening of} $X^{ab}$.

\begin{lem}\label{NC-sm-center} 
(i) Let $R$ be a $d$-smooth algebra, such that $\dim R^{ab}\ge 2$ and $R^{ab}$ is connected.
Assume that $d\ge 1$. Then the center of $R$ is
$\C+\I_{d+1}R$.

\noindent
(ii) Let $\OO^{NC}_X$ be an NC-smooth thickening of a smooth scheme $X$, where $\dim X\ge 2$. 
Then the center of $\OO^{NC}_X$ is the constant sheaf $\C_X$.
\end{lem}

\Pf . (i) Let $Z(R)$ denote the center of $R$.
We have a central extension of algebras
$$0\to \I_{d+1}R\to R\to R'=R/\I_{d+1}R\to 0,$$
hence, we have the inclusion $\C+\I_{d+1}R\sub Z(R)$.
In the case $d=1$ we have $R'=R^{ab}$ and the commutator pairing associated with the above extension
is 
$$[\wt{f},\wt{g}]=df\we dg\in \Om^2_{R^{ab}}\simeq \I_2R,$$ 
where $\wt{f},\wt{g}\in R$ are lifts of $f,g\in R'=R^{ab}$.
This easily implies that
an element of $Z(R)$ projects to $\C\sub R'$.

In the case $d>1$, by the induction assumption, we can assume that
$Z(R')=\C+\I_dR'$. Hence, it is enough to investigate elements of $Z(R)$ that project to elements of $\I_dR'\sub R'$.
Let us consider the commutator pairing
$$\I_dR'\times R^{ab}\to \I_{d+1}R: (\a,f)\mapsto [\wt{\a},\wt{f}],$$
where $\wt{\a},\wt{f}\in R$ are lifts of $\a$ and $f$. 
We claim that this pairing is induced by the natural commutator pairing
$$U(\LLie_+(\Om^1_{R^{ab}}))_d \times \Om^1_{R^{ab}}\to U(\LLie_+(\Om^1_{R^{ab}}))_{d+1},$$
where we use the notation of \cite[Sec.\ 2.1]{PT} (in particular, $\LLie_+(?)$ denotes the degree $\ge 2$ part of the free Lie algebra), and an isomorphism
$$\gr^n_{\I}(R)\simeq U(\LLie_+(\Om^1_{R^{ab}}))_n$$
for $n\le d+1$ (see \cite[Cor.\ 2.3.15]{PT}).
More precisely, we claim that
\begin{equation}\label{commutator-formula-eq}
[\wt{\a},\wt{f}]=-[\a,df]_U,
\end{equation}
where we view $\a$ as an element of $\I_dR'=\gr^d_{\I}(R)$, and 
in the right-hand side we take the commutator in the algebra $U(\LLie_+(\Om^1_{R^{ab}}))$.
Indeed, by \cite[Cor.\ 2.3.9]{PT}, we can realize $R$ as a subalgebra in $T(\Om^1_{R^{ab}})/T^{\ge d+2}(\Om^1_{R^{ab}})$
(where $T(?)$ denotes the tensor algebra over $R^{ab}$),
so that the projection $R\to R^{ab}$ is induced by the projection to $T^0=R^{ab}$.
Furthermore, the elements in the image of $R$ have tensor components of the form $(f,-df,\ldots)$.
Since $T^0$ is in the center of the tensor algebra, this immediately implies formula \eqref{commutator-formula-eq}.
Thus, if $\a\in\I_dR'$ lifts to an element of $Z(R)$ then $[\a,df]_U=0$ for any $f$. Since, $U(\LLie_+(\Om^1_{R^{ab}}))$
is a subalgebra in the tensor algebra $T(\Om^1_{R^{ab}})$, this implies that $\a$ is in the center of 
$T(\Om^1_{R^{ab}})$, hence, $\a=0$ (since $\dim R^{ab}\ge 2$). This implies that $Z(R)=\C+\I_{d+1}R$.

\noindent
(ii) It is enough to check this in the case when $X$ is affine connected, i.e., $X$ is a formal spectrum of an NC-smooth
algebra $R$ such that $R^{ab}$ is connected. 
Now the assertion easily follows from (i).
\ed

By a vector bundle $E$ on a NC-nilpotent scheme $X$ we mean a sheaf of right $\OO$-modules which is
locally free of finite rank. We denote by $E^{ab}$ the induced vector bundle on $X^{ab}$.

\begin{lem}\label{triv-bun-lem} Assume that $X\sub X'$ is a nilpotent extension of affine NC-nilpotent
schemes, i.e., $\OO_X$ is a quotient of $\OO_{X'}$ by a nilpotent ideal.
Let $E'$ be a vector bundle over $X'$,
and $E$ the induced vector bundle over $X$.
Let $\varphi:\OO_X^n\to E$ be a trivialization. Then $\varphi$ extends to a trivialization $\OO_{X'}^n\to E$.
\end{lem}

\Pf . It is enough to consider the case when $0\to \II\to \OO_{X'}\to \OO_X\to 0$ is a central extension. Then
$\II$ is a quasicoherent sheaf over $\OO_{X^{ab}}$, so
$$H^1(X',E'\ot\II)=H^1(X^{ab},E^{ab}\ot \II)=0.$$ 
Thus, the $n$ global sections of $E$ defining the trivialization can be lifted to global sections of $E'$.
It is easy to see (arguing locally) that they give a trivialization of $E'$.
\ed

\subsection{The category of affine almost NC schemes}\label{aNC-sec}



The category $a\NN$ has the same objects as $\NN$,
while the morphisms in $a\NN$ are equivalence classes of homomorphisms $A\to B$, where $f_1,f_2:A\to B$ are
equivalent if there exists $b\in B^*$ such that $f_2=bf_1b^{-1}$.
We denote $a\NN_d\sub a\NN$ the full subcategory of NC-nilpotent algebras of degree $d$.

Given a ring $A$ in $\NN$ and a multiplicative set $\ov{S}\sub A^{ab}$, let $S$ denote the preimage of $\ov{S}$ under the projection $A\to A^{ab}$. Then $S$ satisfies Ore conditions and $S^{-1}A$ is again NC-nilpotent (see \cite[Sec.\ 2.1]{Kapranov}).
For any $B\in\NN$ the composition with the localization morphism $\iota:A\to A[S^{-1}]$ induces an embedding
$\Hom_{\NN}(A[S^{-1}],B)\hra \Hom_{\NN}(A,B)$ with the image consisting of $[f]$ such that $f(S)\sub B^*$.
Since the latter condition is invariant with respect to our equivalence relation on $\Hom_{\NN}(A,B)$, the composition
with $[\iota]$ gives an embedding
$$\Hom_{a\NN}(A[S^{-1}],B)\hra \Hom_{a\NN}(A,B)$$
with the same characterization of the image.

Note also that for $B\in\NN$ an element $b\in B$ is invertible if and only if its image in $B^{ab}$ is invertible.
Thus, a homomorphism $f:A\to B$ factors through $A[S^{-1}]$ if and only if the induced homomorphism
$f^{ab}:A^{ab}\to B^{ab}$ factors through $A^{ab}[\ov{S}^{-1}]$, where $\ov{S}\sub A^{ab}$ is the image of $S$.
It follows that we have a cartesian square of sets
\begin{diagram}
h_{A[S^{-1}]}(B) &\rTo{}& h_A(B)\\
\dTo{}&&\dTo{}\\
h_{A^{ab}[\ov{S}^{-1}]}(B^{ab})&\rTo{} &h_{A^{ab}}(B^{ab})
\end{diagram}

Now let $R$ be an NC-complete algebra and let $T\in R^{ab}$ be a multiplicative subset. Following Kapranov 
\cite[Def.\ (2.1.8)]{Kapranov}, we set
$$R[\![T^{-1}]\!]:=\projlim (R/\I_dR)[T_d^{-1}],$$
where $T_d\sub R/\I_dR$ is the preimage of $T$. In the case when $T=\{f^n \ |\ n\ge 0\}$, for some element $f\in R^{ab}$,
we denote the above algebra simply as $R[\![f^{-1}]\!]$.

For an NC-complete algebra $R$ we denote by $\ov{h}_R$ the corresponding functor on $a\NN$:
$\ov{h}_R(B)$ is the set of conjugacy classes of algebra homomorphisms $R\to B$. 
Since the images of both horizontal arrows in the above cartesian square 
are stable under the action of inner automorphisms of $B$,
we deduce that the similar square
\begin{equation}\label{localization-square}
\begin{diagram}
\ov{h}_{R[\![T^{-1}]\!]}(B) &\rTo{}& \ov{h}_R(B)\\
\dTo{}&&\dTo{}\\
h_{R^{ab}[T^{-1}]}(B^{ab})&\rTo{} &h_{R^{ab}}(B^{ab})
\end{diagram}
\end{equation}
is still Cartiesian for any $B\in\NN$.

Let $a\NN\CC$ denote the category of NC-complete algebras with morphisms given by algebra homomorphisms
viewed up to conjugation, i.e., up to post-composing with an inner automorphism.
We denote by $a\NN\CC\SS_{is}$ the subcategory in $a\NN\CC$, whose objects are {\it NC-smooth} algebras, with {\it isomorphisms} in $a\NN\CC$ as morphisms.

\begin{lem}\label{aNC-Ioneda-lem} 
The functor 
$$a\NN\CC\SS_{is}^{op}\to \Fun_{is}(a\NN,\Sets): R\mapsto \ov{h}_R$$
is fully faithful, where $\Fun_{is}$ is the category of functors and natural isomorphisms between them.
\end{lem}

\Pf . Note that for any $d\ge 0$, the restriction $\ov{h}_R|_{a\NN_d}$ is naturally isomorphic to the representable
functor $\ov{h}_{R/\I_{d+2}R}$. Thus, for NC-complete algebras $R$ and $R'$, we have a natural identification
$$\Isom(\ov{h}_{R'},\ov{h}_{R})\simeq \projlim_d \Isom_{a\NN}(R/\I_{d+2}R,R'/\I_{d+2}R'),$$
where $\Isom_{a\NN}(?,?)$ denotes the set of isomorphisms in the category $a\NN$.
Thus, it suffices to prove that if $R$ and $R'$ are NC-smooth then the natural map
\begin{equation}\label{projlim-eq}
\Isom_{a\NN\CC}(R,R')\to \projlim_d \Isom_{a\NN}(R/\I_{d+2}R,R'/\I_{d+2}R')
\end{equation}
is a bijection. To check surjectivity, assume we are given a collection of algebra homomorphisms 
$$f_d:R/\I_{d+2}R\to R'/\I_{d+2}R',$$
which are compatible up to conjugation, i.e., the homomorphism $f_{d+1,d}:R/\I_{d+2}R\to R'/\I_{d+2}R'$ induced
by $f_{d+1}$ is equal to $\th_{u_d}f_d$, where $\th_{u_d}$ is the inner automorphism associated with
a unit $u_d\in R'/\I_{d+2}R'$. Now, starting from $d=0$, we can recursively correct $f_{d+1}$ by an inner automorphism
of $R'/\I_{d+3}R'$, so that the homomorphisms $(f_d)$ become compatible on the nose (not up to an inner automorphism).
Since $R'$ is NC-complete, this defines a unique homomorphism $f:R\to R'$ inducing $(f_d)$. Furthermore, since
$R$ is NC-complete, we see that $f$ is an isomorphism if and only if all $f_d$ are isomorphisms.

It remains to check that \eqref{projlim-eq} is injective. Thus, given two isomorphisms $f,f':R\to R'$ such that
the induced isomorphisms $f_d$ and $f'_d$ are conjugate for each $d$, we have to check that $f$ and $f'$ are conjugate.
By considering $f^{-1}f'$, we reduce the problem to checking that if we have an automorphism $f:R\to R$ such that
$f_d$ is an inner automorphism of $R/\I_{d+2}R$ for each $d$, then $f$ is inner. For any algebra $A$, 
let us denote by $\Inn(A)$ the group of inner automorphisms of $A$. Note that we have an exact sequence
of groups
$$1\to Z(A)^*\to A^*\to \Inn(A)\to 1.$$
Applying this to each algebra $R/\I_{d+2}R$, and passing to projective limits, we have an exact sequence
$$1\to \projlim_d Z(R/\I_{d+2}R)^*\to \projlim_d (R/\I_{d+2}R)^*\rTo{\rho} \projlim_d \Inn(R/\I_{d+2}R).$$
We claim that the arrow $\rho$ in this sequence is surjective. Indeed, it is enough to check that the
inverse system $(Z(R/\I_{d+2}R)^*)$ satisfies the Mittag-Leffler condition. But by Lemma \ref{NC-sm-center}(i),
for $d\ge 1$, the image of the projection 
$$Z(R/\I_{d+2}R)^*\to Z(R/\I_{d+1}R)^*$$
is equal to $\C^*$, which implies the required stabilization. Thus, the map $\rho$ is surjective. Note that the source of this
map can be identified with $R^*$. Thus, we deduce the surjectivity of the natural map
$$R^*\to \projlim_d \Inn(R/\I_{d+2}R).$$
Hence, we can compose $f$ with an inner automorphism $\th_u$ of $R$, such that $f'=\th_u f$ induces the identity
automorphism of $R/\I_{d+2}R$ for each $d$. It follows that $f'=\id$, i.e., $f$ is inner.
\ed


\subsection{Gluing}\label{gluing-sec}

We can define Zariski topology on $a\NN^{op}$ naturally. However, this is not a subcanonical topology, i.e.,
representable functors are not necessarily sheaves with respect to this topology.
Namely, suppose $f_1,f_2:A\to B$ is a pair of homomorphisms, inducing the same homomorphism $A^{ab}\to B^{ab}$.
Assume also that we have a covering of $\Spec(A^{ab})$ by distinguished affine opens,
$\Spec(A_{\ov{g}_i})$ such that $f_1$ and $f_2$ become conjugate as morphisms from $A_{g_i}$ to $B_{g_i}$.
It may happen that $f_1$ and $f_2$ are still not conjugate by an element of $B^*$. 

\begin{ex}
Let $R:=k[x,y]$. For any ideal $I\sub R$ we can consider the central extension $A$ of $R$ by $\Om^2_{R/k}\ot_R R/I$,
obtained from the universal central extension via the natural homomorphism $\Om^2_{R/k}\to \Om^2_{R/k}\ot_R R/I$.
We consider a pair of homomorphisms 
$$f_1=\id, f_2=\id+\de: A\to A,$$
where $\de:A\to R\to \Om^2_{R/k}\ot_R R/I$ is a derivation given by
$$\de(r)=\om\we dr \mod I\Om^2_{R/k},$$
for some $1$-form modulo $I$, $\om\in \Om^1_{R/k}\ot_R R/I$. 
We are going to prove that for $I=(xy-1)$, there exists a $1$-form $\om$ such that $f_1$ and $f_2$ are locally conjugate, but
not globally conjugate. Note that if $R_g$ is a localization of $R$ then the corresponding localization of $A$
is a central extension of $R_g$ by $\Om^2_{R_g/k}\ot_{R_g} (R/I)_g$.
It is easy to see that the condition for $f_1$ and $f_2$ to be conjugate over $\Spec(R_g)$ is that for some $r\in R_g^*$ one has
$$\de(r')=r^{-1}dr\wedge dr' \mod I\Om^2_{R_g/k}$$
for any $r'\in R$. Since the morphism $\eta\mapsto \eta\wedge ?$ gives an isomorphism
$$\Om^1_{R_g/k}\ot_{R_g} (R/I)_g\simeq \Hom_{R_g}(\Om^1_{R_g/k},\Om^2_{R_g/k}\ot_{R_g} (R/I)_g),$$
this is equivalentl to the condition 
$$\om\equiv r^{-1}dr \mod I\Om^1_{R_g/k}.$$
Let us consider the homomorphism of sheaves on $\Spec(R)=\A^2_k,$
$$\tau:\OO^*\to \Om^1\ot_{\OO} \OO/I,$$
induced by $\phi\mapsto \phi^{-1}d\phi$.
Thus, the condition on $\om$ means that it comes from a global section of the sheaf image $\im(\tau)$, 
but is not in the image of the induced morphism on global sections. 
Since $H^0(\A^2,\OO^*)=k$, the latter condition
is equivalent to $\om\neq 0$.
Now we observe that the sheaf $\Om^1\ot_{\OO} \OO/I$ is supported on the curve $xy=1$ which is contained
in the affine open subset $x\neq 0$. Hence, the $1$-form $dx/x$ gives a well defined nonzero global section of
$\im(\tau)$, as required. 
\end{ex}

Because of this we do not try to glue affine almost NC schemes using sheaves on $a\NN$.
Instead, we show that a locally representable formally smooth functor on $a\NN$ always leads 
to an algebroid over the underlying commutative smooth scheme $X$, that corresponds locally to an NC-smooth 
thickening of $X$.

Recall that a {\it $\C$-algebroid} $\AA$ 
over a topological space $X$ is a stack of $\C$-linear categories over $X$, such that $\AA$ is locally non-empty and 
any two objects of $\AA(U)$ are locally isomorphic. We refer to \cite[Sec.\ 2.1]{KS} for basic results on algebroids.

\begin{defi}\label{NC-stack-alg-defi} 
Let $X$ be a smooth scheme. An {\it NC-smooth algebroid thickening} of $X$ is a $\C$-algebroid 
$\AA$ over $X$
such that for every object $\si\in \AA(U)$ over an open subset $U\sub X$ the sheaf of algebras 
$\Eend_\AA(\si)$ is an NC-smooth thickening of $U$.
\end{defi}

For a sheaf of $\C$-algebras $A$ over $X$ we have the corresponding $\C$-algebroid with a fixed global object
$\si$ such that $A$ is the endomorphism algebra of $\si$.

\begin{defi}
For a $\C$-algebroid $\AA$, we define the {\it center} of $\AA$ as the sheaf 
$$\ZZ_\AA:=\Eend(\Id_\AA)$$
of endomorphisms of the identity functor on $\AA$. 
We say that a $\C$-algebroid $\AA$ has {\it trivial center} if the natural map
 of sheaves $\C_X\to \ZZ_\AA$ is an isomorphism. 
\end{defi}

It is easy to see that for any local object $\si\in\AA(U)$ 
one has a natural identification of $\ZZ_\AA|_U$ with the center of the sheaf of algebras $\Eend_\AA(\si)$.
Thus, by Lemma \ref{NC-sm-center}(ii), any NC-smooth algebroid thickening has trivial center.

We are going to prove a general gluing result for sheaves of $\C$-algebras with trivial centers and then apply them to
construct NC-smooth algebroid thickenings.

\begin{lem}\label{C*-algebroid-lem}
(i)  Let $\AA$ and $\AA'$ be a pair of $\C$-algebroids with trivial centers over an irreducible scheme $X$, and 
let $F,G:\AA\to \AA'$ be a pair of equivalences.
Assume that for an open covering $(U_i)$ of $X$ we have an isomorphism
$F|_{U_i}\simeq G|_{U_i}$. Then there exists an isomorphism $F\simeq G$.

\noindent
(ii) Let $\AA$ and $\AA'$ be a pair of $\C$-algebroids with trivial centers over an irreducible scheme $X$.
Assume that for an open covering $(U_i)$ of $X$ we have an equivalence 
$$F_i:\AA|_{U_i}\to \AA'|_{U_i}$$
and that for each pair $i,j$, we have an isomorphism 
$$F_i|_{U_{ij}}\simeq F_j|_{U_{ij}},$$
where $U_{ij}=U_i\cap U_j$.
Then there exists an equivalence $F:\AA\to \AA'$ such that $F|_{U_i}\simeq F_i$.
Such an equivalence is unique up to an isomorphism.

\noindent
(iii) Let $U_i$ be an open covering of an irreducible scheme $X$, and for each $i$ let $\AA_i$ be a $\C$-algebroid
with trivial center over $U_i$. Assume that for every $i,j$, we have an equivalence 
$$F_{ij}:\AA_i|_{U_{ij}}\to \AA_j|_{U_{ij}},$$ 
such that for every $i,j,k$, there is an isomorphism 
$$F_{jk}|_{U_{ijk}}\circ F_{ij}|_{U_{ijk}}\simeq F_{ik}|_{U_{ijk}},$$
where $U_{ijk}=U_i\cap U_j\cap U_k$. Then there exists a $\C$-algebroid $\AA$ over $X$
and equivalences $F_i:\AA|_{U_i}\to \AA_i$, such that for every $i,j$, there is an isomorphism
$$F_{ij}\circ F_i|_{U_{ij}}\simeq F_j|_{U_{ij}}.$$
Furthermore, such $\AA$ is unique up to an equivalence.
\end{lem}

\Pf . 
(i) Let us choose for each $i$ an isomorphism $\phi_i:F|_{U_i}\to G|_{U_i}$. Then for each $i,j$, we have
$$\phi_j|_{U_{ij}}=\phi_i|_{U_{ij}}\circ c_{ij},$$
where $c_{ij}$ is an autoequivalence of $F_i|_{U_{ij}}$.
Since $F_i$ is an equivalence, we have $\und{\Aut}(F)\simeq \und{\Aut}(\id_{\AA})$.
Locally, the sheaf $\und{\Aut}(\id_{\AA})$ is given by the center of $\End_\AA(\si)$, where $\si$ is an object of $\AA$.
Hence, by Lemma \ref{NC-sm-center},
the natural morphism of sheaves $\C^*_X\to \und{\Aut}(\id_\AA)$ is an isomorphism.
Thus, $c_{ij}$ is a Cech $1$-cocycle with values in $\C^*_X$. Since $X$ is irreducible, the corresponding Cech cohomology
is trivial, so we can multiply $\phi_i$ by appropriate constants in $\C^*$, to make them compatible on double intersections.
The corrected isomorphisms glue into a global isomorphism $F\to G$.

\noindent
(ii) Let us choose for each $i,j$ an isomorphism $\phi_{ij}:F_i|_{U_{ij}}\to F_j|_{U_{ij}}$.
Then for each $i,j,k$, the composition $c_{ijk}=\phi_{ki}\phi_{jk}\phi_{ij}$ is an autoequivalence of $F_i|_{U_{ijk}}$,
where $c_{ijk}$ is a Cech $2$-cocycle with values in $\C^*_X$. As above, choosing representation of $c_{ijk}$ as a coboundary
allows to correct $\phi_{ij}$ by constants in $\C^*$, so that the isomorphisms $\phi_{ij}$ are compatible on triple intersections.
Hence, we can glue $(F_i)$ into the required global equivalence $F:\AA\to\AA'$.
The fact that $F$ is unique up to an isomorphism follows from (i).

\noindent
(iii) For every $i,j,k$, let us choose an isomorphism
$$g_{ijk}:F_{jk}|_{U_{ijk}}\circ F_{ij}|_{U_{ijk}}\to F_{ik}|_{U_{ijk}}.$$
Then for every $i,j,k,l$, we have over $U_{ijkl}$, 
$$g_{ikl} (F_{kl}*g_{ijk})=c_{ijkl} g_{ijl} (g_{jkl}*F_{ij})$$
for some $c_{ijkl}\in\und\Aut(F_{il})(U_{ijkl})=\C^*$. Furthermore, $(c_{ijkl})$ is a Cech $3$-cocycle with values in $\C^*_X$.
Hence, we can mulitply $g_{ijk}$ with appropriate constants to make them compatible on quadruple intersections.
This allows to glue $(\AA_i)$ into a global $\C$-algebroid $\AA$ over $X$ (see \cite[Prop.\ 2.1.13]{KS}).
The uniqueness of $\AA$ up to an equivalence follows from (ii).
\ed

\begin{prop}\label{general-gluing-prop} 
Let $(U_i)$ be an open covering of an irreducible scheme $X$. Assume that for each $i$ we are given
a sheaf of $\C$-algebras $A_i$ with trival center over $U_i$, and for each pair $i<j$, a covering $(V_k=V_{ij,k})$ of $U_{ij}$,
together with isomorphisms of sheaves of $\C$-algebras
$$\a_{ij,V_k}:A_i|_{V_k}\to A_j|_{V_k}$$
for all $k$.
We assume that the restrictions of $\a_{ij,V_k}$ and $\a_{ij,V_l}$ to $V_k\cap V_l$ differ by an inner automorphism.
Also, we assume that for $i<j<k$ there exists a covering $(W_l=W_{ijk,l})$ of $U_{ijk}$ such that
$\a_{jk}|_{W_l}\circ \a_{ij}|_{W_l}$ and $\a_{ik}|_{W_l}$ differ by an inner automorphism.
Then there exists a $\C$-algebroid $\AA$ over $X$, together with equivalences of $\C$-algebroids, 
$$F_i:\AA|_{U_i}\to \AA_i,$$
where $\AA_i$ is the $\C$-algebroid over $U_i$ associated with $A_i$,
such that for $i<j$ there exist isomorphisms 
$$\a_{ij,V_k}\circ F_i|_{V_k}\simeq F_j|_{V_k}$$ 
over the covering $(V_k)$ of $U_{ij}$. Such an algebroid $\AA$ is unique up to an equivalence.
\end{prop}

\Pf . Each isomorphism $\a_{ij,V_k}$ gives an equivalence 
$$F_{ij,V_k}:\AA_i|_{V_k}\to \AA_j|_{V_k}.$$ 
Since the local autoequivalence of $\AA_i$ associated with an inner automorphism of $A_i$ is isomorphic to the identity, we get that $F_{ij,V_k}$ and $F_{ij,V_l}$ induce isomorphic equivalences over $V_k\cap V_l$.
By Lemma \ref{C*-algebroid-lem}(ii), we obtain an equivalence defined over $U_{ij}$,
$$F_{ij}:\AA_i|_{U_{ij}}\to \AA_j|_{U_{ij}},$$
such that $F_{ij}|_{V_k}\simeq\a_{ij,V_k}$.

Furthermore, we claim that over $U_{ijk}$ there is an isomorphism 
\begin{equation}\label{F-ijk-eq}
F_{jk}|_{U_{ijk}}\circ F_{ij}|_{U_{ijk}}\simeq F_{ik}|_{U_{ijk}}.
\end{equation}
Indeed, by assumption, we have a similar isomorphism over each open subset from the covering $(W_l)$ of $U_{ijk}$.
Thus, our claim follows from Lemma \ref{C*-algebroid-lem}(i), applied to the equivalences in
both sides of \eqref{F-ijk-eq}.

Finally, we can apply Lemma \ref{C*-algebroid-lem}(iii) to conclude the existence and uniqueness
of the required NC-smooth algebroid
$\AA$ over $X$.
\ed

Now we are going to apply the above general result to NC-smooth thickenings.

For a functor $h$ on $a\NN$ such that $h_{\Com}=h_X$ and an open subset $U\sub X$ we define the
subfunctor $h_{/U}\sub h$ by 
$$h_{/U}(\La)=h(\La)\times_{h_{X}(\La^{ab})} h_{U}(\La^{ab}),$$
where we use the identification $h(\La^{ab})\simeq h_X(\La^{ab})$.

\begin{lem}\label{dist-af-lem} 
Let $h=\ov{h}_R$, where $R$ is an NC-complete algebra. Then for any distinguished affine $D(f)\sub \Spec(A^{ab})$
we have an equality of subfunctors $h_{/D(f)}=\ov{h}_{A[\![f^{-1}]\!]}$.
\end{lem}

\Pf . This follows immediately from the cartesian square \eqref{localization-square} with $T=\{f^n \ |\ n\ge 0\}$.
\ed


\begin{lem}\label{gluing-lem} 
Let $h$ be a functor on $a\NN$ such that $h|_{\Com}=h_X$ for some scheme $X$.
Assume that $(U_i)$ is an affine covering of $X$, such that for every $i$ we have
an isomorphism $h_{/U_i}\simeq \ov{h}_{A_i}$ for some 
$A_i\in \NN$. Let us denote also by $A_i$ the corresponding sheaf of algebras over $U_i$.
Then for every open subset $V\sub U_i\cap U_j$, which is distingushed in both $U_i$ and $U_j$, we
have an isomorphism 
$$\a_{ij,V}:A_i|_V\simeq A_j|_V$$ 
compatible with the isomorphisms
$\ov{h}_{A_i(V)}\simeq h_{/V}\simeq \ov{h}_{A_j(V)}$. Furthermore, for another such open $V'\sub U_i\cap U_j$
the isomorphisms $\a_{ij,V}|_{V\cap V'}$ and $\a_{ij,V'}|_{V\cap V'}$ differ by an inner automorphism.
Also, for any open $V\sub U_i\cap U_j\cap U_k$, distinguished in $U_i$, $U_j$ and $U_k$,
we have 
$$\a_{jk}|_V\circ \a_{ij}|_V=\a_{ik}|_V\circ \Ad(u_{ijk})$$ 
for some $u_{ijk}\in A_i(V)^*$. 
\end{lem}
 
\Pf . Let us fix an isomorphism $h_{/U_i}\simeq\ov{h}_{A_i}$ for each $i$.
Suppose $V\sub U_i\cap U_j$ is a distinguished affine open in both $U_i$ and $U_j$. Then 
$$\ov{h}_{A_i,/V}\simeq h_{/V}\simeq \ov{h}_{A_j,V}.$$
Thus, by Lemmas \ref{dist-af-lem} and \ref{aNC-Ioneda-lem}, 
we have an isomorphism between the corresponding localizations of $A_i$ and $A_j$ in $a\NN$,
and hence, an isomorphism $\a_{ij}:A_i|_V\simeq A_j|_V$, defined uniquely up to an inner automorphism.
For $V\sub U_i\cap U_j\cap U_k$ the compatibility between $\a_{ij}$, $\a_{jk}$ and $\a_{ik}$, up to an inner automorphism,
follows from the compatibility of all of these isomorphisms with the isomorphisms of $\ov{h}_{A_i,/V}$, $\ov{h}_{A_j,/V}$ and
$\ov{h}_{A_k,/V}$ with $h_{/V}$.
\ed



\begin{thm}\label{functor-algebroid-thm} 
Let $h$ be a formally smooth
functor on $a\NN$ such that $h|_{\Com}=h_X$ and $h$ is locally representable, i.e., there exists an open affine covering
$(U_i)$ of $X$, and isomorphisms 
$$h_{/U_i}\simeq \ov{h}_{A_i},$$ 
where $A_i$ is an NC-smooth thickening of $U_i$.
Then there exists an NC-smooth algebroid $\AA$ over $X$ and equivalences of algebroids
$$F_i:\AA|_{U_i}\to A_i,$$ 
such that for every open subset $V\sub U_i\cap U_j$, distinguished in both $U_i$ and $U_j$, 
there is an isomorphism
$$g_{ij}\circ F_i|_V\simeq F_j|_V,$$
where $g_{ij}:A_i|_V\to A_j|_V$ is a representative (up to conjugation) of the isomorphism
$\ov{h}_{A_i|_V}\simeq h_{/V}\simeq \ov{h}_{A_j|_V}$.
\end{thm}

\Pf . Without loss of generality we can assume that $X$ is connected.

First, we apply Lemma \ref{gluing-lem} and obtain isomorphisms 
$$\a_{ij,V}:A_i|_V\to A_j|_V$$ 
for every open $V\sub U_i\cap U_j$, distinguished in both $U_i$ and $U_j$, such that these isomorphisms for $V$ and $V'$
and for $V\sub U_i\cap U_j\cap U_k$, are compatible up to an inner automorphism.
Hence, we are in the setup of Proposition \ref{general-gluing-prop}, where as open coverings of $U_{ij}$ (resp., $U_{ijk}$)
we take the covering by all open affines which are distinguished in $U_i$ and $U_j$ (resp., $U_i$, $U_j$ and $U_k$).
Note that the centers of $A_i$ are trivial by Lemma \ref{NC-sm-center}(ii).
Thus, applying Proposition \ref{general-gluing-prop} we get the required NC-smooth algebroid over $X$.
\ed

\subsection{Recollections on nonabelian $H^1$}\label{nonab-H1-sec}

We are going to use some basic constructions involving nonabelian cohomology, which we recall here.
The comprehensive reference is in Giraud's book \cite{Gir} (more specifically, we use \cite[Sec.\ 3.3,3.4]{Gir}). 
A more explicit treatment in terms of Cech cocycles
is given in \cite[Sec.\ 2.6.8]{Manin}, however, it contains one mistake that we will correct.

For a sheaf of groups $\GG$ on a topological space $X$ and an open covering $\UU=(U_i)$ of $X$, the set of $1$-cocycles
 $Z^1(\UU,\GG)$ consists of $g_{ij}\in\GG(U_{ij})$, such that $g_{ii}=1$, $g_{ij}g_{ji}=1$ and 
$$g_{ij}|_{U_{ijk}}g_{jk}|_{U_{ijk}}=g_{ik}|_{U_{ijk}}.$$
Two such $1$-cocycles $(g_{ij})$ and $(\wt{g}_{ij})$ are cohomologous if
$$\wt{g}_{ij}=h_i|_{U_{ij}} g_{ij}h_j^{-1}|_{U_{ij}},$$
for some $h_i\in\GG(U_i)$. We denote by $H^1(\UU,\GG)$ the corresponding set of equivalence classes (pointed
by the class of the trivial cocycle). The nonabelian
cohomology set $H^1(X,\GG)$ is obtained by taking the limit over all open coverings.
Note that our convention for 
nonabelian $1$-cocycles is the same as in \cite{Gir} and differs from that of \cite[Sec.\ 2.6.8]{Manin} by passing to inverses.
For brevity, from now on, we stop writing explicitly the restrictions to the intersections in formulas involving sections defined over
different open subsets. 

For a homomorphism $\GG_1\to\GG_2$ the induced map of pointed sets $H^1(X,\GG_1)\to H^1(X,\GG_2)$ is defined
in an obvious way. 
Now assume we are given an abelian extension of sheaves of groups
$$1\to \AA\to \GG'\rTo{p} \GG\to 1.$$
This means that $\AA$ is a sheaf of abelian groups, which is a normal subsheaf in $\GG'$, and $\GG$ is the corresponding
quotient. Then we have a natural connecting map 
$$\de_0:H^0(X,\GG)\to H^1(X,\AA)$$
such that $\de_0(g)=1$ if and only if $g$ lifts to a global section of $\GG'$. Namely,
for an open covering $U_i$ we can find $g'_i\in \GG'(U_i)$, such that $p(g'_i)=g$ and set $\de_0(g)$ to be the class
of the $1$-cocycle 
$(g'_i)^{-1}g'_j\in\AA(U_{ij})$.
Note that $\de_0$ is {\it not a homomorphism} in general. Rather, it satisfies
\begin{equation}\label{de0-crossed-eq}
\de_0(g_1g_2)=g_2^{-1}(\de_0(g_1))+\de_0(g_2),
\end{equation}
where we write the group structure in $H^1(X,\AA)$ additively and
use the natural action of $H^0(X,\GG)$ on $H^1(X,\AA)$ induced by the adjoint action of $\GG$ on $\AA$.
(This means that $g\mapsto \de_0(g^{-1})$ is a {\it crossed homomorphism}.)
An equivalent restatement of \eqref{de0-crossed-eq} is that there is a {\it twisted action} of $H^0(X,\GG)$ on
$H^1(X,\AA)$ given by
\begin{equation}\label{twisted-H0-action-eq}
g \times a= g(a)+\de_0(g^{-1}), \ \text{ where } g\in H^0(X,\GG), a\in H^1(X,\AA).
\end{equation}
Explicitly, the usual action of $g\in H^0(X,\GG)$ on a class of a Cech $1$-cocycle $(a_{ij})$ with values in $\AA$
is given by $g'_ia_{ij}(g'_i)^{-1}$, where $g'_i\in\GG'(U_i)$ are liftings of $g$. On the other hand, the twisted action of
$g$ on $a_{ij}$ is given by $g'_ia_{ij}(g'_j)^{-1}$.

Next, starting from a class $g\in H^1(X,\GG)$ we can construct a class
$$\de_1(g)\in H^2(X,\AA^g)$$
such that $\de_1(g)=0$ if and only if $g$ is in the image of the map
$H^1(X,\GG')\to H^1(X,\GG)$. Here $\AA^g$ is the sheaf obtained from $\AA$ by twisting with $g$.
Namely, if $g$ is represented by a Cech $1$-cocycle $g_{ij}\in \GG(U_i)$ then we have isomorphisms
$\psi_i:\AA|_{U_i}\to \AA^g|_{U_i}$ such that $\psi_j=\psi_i\circ g_{ij}$ over $U_{ij}$.
To construct $\de_1(g)$, for some covering $(U_i)$, we can choose liftings $g'_{ij}\in \GG'(U_{ij})$
for a $1$-cocycle $(g_{ij})$ representing $g$ (such that $g'_{ij}g'_{ji}=1$ and $g'_{ii}=1$). 
Then $\de_1(g)$ is the class of the $2$-cocycle $(\psi_i(g'_{ij}g'_{jk}g'_{ki}))$ with values in $\AA^g$.

Finally, for a given class $g\in H^1(X,\GG)$ we need the following description of the fiber of the map
$$H^1(X,\GG')\rTo{H^1(p)} H^1(X,\GG)$$ 
over $g$. Assume that this fiber is nonempty and let us choose an element $g'\in H^1(X,\GG')$
projecting to $g$. Then we have an exact sequence of twisted groups
$$1\to \AA^g\to (\GG')^{g'}\to \GG^g\to 1.$$
Thus, as before we have two actions of the group $H^0(X,\GG^g)$ on $H^1(X,\AA^g)$. Now we can construct a surjective map
\begin{equation}\label{H1-fibers-map-eq}
H^1(X,\AA^g)\to H^1(p)^{-1}(g),
\end{equation}
such that the fibers of this map are the orbits of the twisted action of $H^0(X,\GG^g)$ on $H^1(X,\AA^g)$ (see
\eqref{twisted-H0-action-eq}).
Namely, let $(g'_{ij})$ be a Cech $1$-cocycle representing $g'$, and let 
$a_{ij}\in \AA(U_{ij})$ be the $g$-twisted $1$-cocycle, so that $\psi_i(a_{ij})$ is a $1$-cocycle with values in $\AA^g$.
This means that over $U_{ijk}$ one has
$$a_{ij} \Ad(g_{ij})(a_{jk})=a_{ik}.$$
Then our map \eqref{H1-fibers-map-eq} sends $(a_{ij})$ to the class of $(a_{ij}g'_{ij})$.

In the particular case when the (usual) action of $H^0(X,\GG^g)$ on $H^1(X,\AA^g)$ is trivial, the corresponding
connecting map 
$$\de_0:H^0(X,\GG^g)\to H^1(X,\AA^g)$$
is a group homomorphism, and the map \eqref{H1-fibers-map-eq} induces an identification of
the cokernel of this homomorphism with $H^1(p)^{-1}(g)$. Equivalently, in this case the map \eqref{H1-fibers-map-eq}
corresponds to a transitive action of $H^1(X,\AA^g)$ on $H^1(p)^{-1}(g)$, such that the stabilizer of any element
is the image of $\de_0$. (In \cite[Sec.\ 2.6.8]{Manin} it is stated incorrectly that such an action exists in the general case.)

\subsection{The functor of NC-families extending a given excellent family}\label{versal-sec}


Let $Z$ be a projective algebraic variety, $B$ a smooth algebraic variety, and let $\EE$ be an excellent family of bundles on
$Z$ with the base $B$ (see Definition \ref{versal-fam-defi}.
Note that our definition is slightly stronger than \cite[Def.\ (5.4.1)]{Kapranov} in that we add condition (d), which is used crucially
in the base change calculations.

For an NC-nilpotent scheme $X$ and a usual scheme $Z$ there is a natural product operation which gives 
an NC-nilpotent scheme $X\times Z$, so that functions on $Z$ become central in $\OO_{X\times Z}$.
In the affine case this corresponds to the operation of extension of scalars $R\mapsto R\ot_\C S$
from NC-nilpotent $\C$-algebras
to NC-nilpotent $S$-algebras, where $S$ is a commutative $\C$-algebra.

Following \cite{Kapranov} we consider the following functor of noncommutative families of vector bundles compatible
with $\EE$.


\begin{defi}\label{fun-defi} 
For an excellent family $\EE$ over a smooth (commutative) base $B$, 
we define the functor $h^{NC}_B: \NN \to Sets$ sending $\Lambda \in \NN$ to the isomorphism classes of objects in the following category $\CC_\Lambda$. Consider
NC-schemes $X=\Spec(\La)$ and $X\times Z$. Let us denote by $X^{ab}_0=\Spec(\La^{ab}_0)$
the reduced scheme associated with
the abelianization of $X$. Then the objects of $\CC_\La$ are the triples $(f, E_\La, \phi)$ consisting of
	\begin{enumerate}
		\item[(i)] a morphism $f: X^{ab}_0 \to B$ of schemes,
		\item[(ii)] a locally free sheaf of right $\OO_{X\times Z}$-modules $E_\La$,
		\item[(iii)] an isomorphism $\phi: \OO_{X^{ab}_0\times Z}\otimes E_\La \overset{\sim}\to (f\times \id)^*\EE$.
	\end{enumerate}
A morphism $(f_1,E_1, \phi_1) \to (f_2, E_2, \phi_2)$ exists only if $f_1 = f_2$ and is given by an isomorphism $E_1\to E_2$
 commuting with the $\phi_i$. On morphisms $h^{NC}_B$ is the usual pullback.
\end{defi}

The following result is stated in \cite{Kapranov} (see \cite[Prop.\ (5.4.3)(a)(b)]{Kapranov}).
However, we believe our stronger assumptions on the family $\EE$, including condition (d), are needed for it to hold,
and we will give a complete proof below. 

\begin{prop}\label{h-smooth-prop}
The functor $h^{NC}_B$ is formally smooth and the natural morphism of functors
$h_B\to h^{NC}_B|_{\Com}$ is an isomorphism.
\end{prop}

\begin{lem}\label{comm-End-lem} 
For any commutative algebra $\La$ and any $(f,E_\La,\phi)\in h^{NC}_B(\La)$ the natural
map 
$$\La\to \End(E_\La)$$
is an isomorphism.
\end{lem}

\Pf . We prove this by the degree of nilpotency of the nilradical of $\La$.
Assume first that $\La$ is reduced. Then we have $E_\La=(f\times\id)^*\EE$.
Hence, by the base change theorem,
\begin{align*}
&H^0(X\times Z,(f\times\id)^*\Eend(\EE))\simeq H^0(X, Rp_{X,*}(f\times\id)^*\Eend(\EE))\simeq\\
&H^0(X, \HH^0(Lf^*R\rho_*\Eend(\EE))),
\end{align*}
where $X=\Spec(\La)$.
Since $R^i\rho_*\Eend(\EE)$ are locally free for $i\ge 1$, we have
$$\HH^0(Lf^*R\rho_*\Eend(\EE))\simeq f^*\rho_*\Eend(\EE)\simeq \OO_X,$$
where in the last isomorphism we used assumption (a). This shows that our assertion holds for such $\La$.

Next, assume we have a central extension $0\to I\to \La'\to \La\to 0$ of commutative algebras,
such that $I$ is a module over $\La_0$, the quotient of $\La$ by its nilradical.
Assume that $\La\to \End(E_\La)$ is an isomorphism for any $(f,E_\La,\phi)\in h^{NC}_B(\La)$
and let us prove a similar statement over $\La'$. Given $(f,E_{\La'},\phi')\in h^{NC}_B(\La')$,
let $E_\La$ be the induced locally free sheaf over $\Spec(\La)\times Z$.
Then we have an exact sequence of coherent sheaves on $\Spec(\La')\times Z$,
$$0\to \EE_{\La_0}\ot p_1^*\II\to \EE_{\La'}\to \EE_\La\to 0,$$
where $\II$ is the ideal sheaf on $\Spec(\La')$ corresponding to $I$.
Taking sheaves of homomorphisms  from $\EE_{\La'}$ we get an exact sequence
$$0\to \Eend(\EE_{\La_0})\ot p_1^*\II\to \Eend(\EE_{\La'})\to \Eend(\EE_\La)\to 0$$
Passing to global sections we obtain a morphism of exact sequences
\begin{diagram}
0 &\rTo{}& I &\rTo{}&\La'&\rTo{}  &\La&\rTo{} 0\\
&&\dTo{}&&\dTo{}&&\dTo{}\\
0\ &\rTo{}& H^0(X^0\times Z,\Eend(\EE_{\La_0})\ot p_1^*\II)&\rTo{}&\End(E_{\La'})&\rTo{}&\End(E_\La)
\end{diagram}
Note that $\EE_{\La_0}\simeq (f\times\id)^*\EE$, so as before we get
\begin{align*}
&H^0(X^0\times Z,\Eend(\EE_{\La_0})\ot p_1^*\II)\simeq H^0(X^0, \II\ot \HH^0(Lf^*R\rho_*\Eend(\EE)))\simeq\\
&H^0(X^0, \II\ot f^*\rho_*\Eend(\EE))\simeq I,
\end{align*}
where $X^0=\Spec(\La_0)$.
Thus, in the above morphism of exact sequences the leftmost and the rightmost vertical arrows are isomorphisms.
Hence, the middle vertical arrow is also an isomorphism.
\ed

\noindent
{\it Proof of Proposition \ref{h-smooth-prop}.} Assume we are given a central extension 
\begin{equation}\label{central-ext-0-eq}
0\to I\to \La'\to \La\to 0 
\end{equation}
in $\NN$ and an element $(f,E_\La,\phi)\in h^{NC}_B(\La)$,
so that $E_\La$ is a locally free sheaf of right $\OO_{X\times Z}$-modules of rank $r$, where $X=\Spec(\La)$. We have to check that it lifts to a locally free sheaf of right $\OO_{X'\times Z}$-modules, where $X'=\Spec(\La')$.
Furthermore, it is enough to consider central extensions as above, where the nilradical of $\La^{ab}$ acts trivially on $I$,
so that $I$ is a $\La^{ab}_0$-module.

We have a natural abelian extension of sheaves of groups on $X^{ab}\times Z$,
\begin{equation}\label{GL-ex-seq}
1\to \Mat_r(\OO_{X^{ab}\times Z})\otimes p_1^*\II\to \GL_r(\OO_{X'\times Z})\to \GL_r(\OO_{X\times Z})\to 1
\end{equation}
where $\II$ is the coherent sheaf on $X^{ab}$ corrresponding to $I$.
The isomorphism class of $E_\La$ corresponds to an element of the nonabelian cohomology 
$H^1(X^{ab}\times Z, \GL_r(\OO_{X\times Z}))$.
By the standard formalism (see Sec.\ \ref{nonab-H1-sec})
the obstruction to lifting this class to a class in $H^1(X^{ab}\times Z, \GL_r(\OO_{X'\times Z}))$
lies in $H^2(X^{ab}\times Z, \Eend(E_{\La^{ab}_0})\otimes p_1^*\II)$, where $E_{\La^{ab}_0}$ is induced by $E_\La$. 
We claim that this group $H^2$ vanishes.
Indeed, we have $E_{\La^{ab}_0}\simeq (f\times \id)^*\EE$.
Applying the base change theorem we get an isomorphism
$$R\Ga(X_0^{ab}\times Z,(f\times\id)^*\Eend(\EE)\ot p_1^*\II)\simeq
R\Ga(X_0^{ab}, \II\ot Lf^*R\rho_*\Eend(\EE)).$$
It remains to observe that by our assumptions (c) and (d), the complex of sheaves
$Lf^*R\rho_*\Eend(\EE)$ has no cohomology in degrees $\ge 2$.

To prove the second assertion we argue by induction on the degree of nilpotency of the nilradical of a test algebra $\La$.
Thus, we consider a square zero extension \eqref{central-ext-0-eq} of commutative algebras, where $I$ is a $\La^{ab}_0$-module, 
and study the corresponding commutative square
\begin{equation}\label{h-B-La'-La-square}
\begin{diagram}
h_B(\La') &\rTo{}& h_B(\La)\\
\dTo{}&&\dTo{}\\
h^{NC}_B(\La') &\rTo{} & h^{NC}_B(\La)
\end{diagram}
\end{equation}
We assume that the right vertical arrow is an isomorphism and we would like to prove the same about the left vertical arrow.
We know that both horizontal arrows are surjective. Furthermore, using the interpretation in terms of nonabelian $H^1$
and the exact sequence \eqref{GL-ex-seq} we can get a description of the preimage of an element $E_\La\in h^{NC}_B(\La)$
under the bottom arrow. Namely, the corresponding sequence of twisted sheaves is
\begin{equation}\label{Aut-ex-seq}
0\to \Eend(E^{ab})\otimes p_1^*\II\to \und{\Aut}(E_{\La'})\to \und{\Aut}(E_\La)\to 1.
\end{equation}
By Lemma \ref{comm-End-lem}, we have $\Aut(E_\La)=\La^*$, and it is easy to see that this group acts trivially on
$H^1(X^{ab}\times Z, \Eend(E_{\La^{ab}})\otimes p_1^*\II))$ (since $\La'$ is in the center of $\und{\Aut}(E_{\La'})$).
It follows that the preimage of $E_\La$ in $h^{NC}_B(\La')$
 is the principal homogeneous space for the abelian group
$$\coker(\Aut(E_\La)\rTo{\de_0} H^1(X^{ab}\times Z, \Eend(E_{\La^{ab}})\otimes p_1^*\II)),$$
where $\de_0$ is the connecting homomorphism associated with \eqref{Aut-ex-seq}. 
However, by Lemma \ref{comm-End-lem}, fixing a lifting $E_{\La'}\in h^{NC}_B(\La')$,
we get that the previous map in the long exact sequence, $\Aut(E_{\La'})\to \Aut(E_\La)$ is just the projection
$(\La')^*\to \La^*$, so it is surjective. This implies that the preimage of $E_\La$ is the principal homogeneous space for
$$H^1(X^{ab}\times Z, \Eend(E_{\La^{ab}_0})\otimes p_1^*\II)\simeq H^0(X_0^{ab}, \II\ot \HH^1(Lf^*R\rho_*\Eend(\EE))).$$
By our assumptions (c) and (d), we have
$$\HH^1(Lf^*R\rho_*\Eend(\EE))\simeq f^*R^1\rho_*\Eend(\EE),$$
thus, the above group is $H^0(X_0^{ab}, I\ot f^*R^1\rho_*\Eend(\EE))$.

On the other hand, different extensions of $\Spec(\La)\to B$ to $\Spec(\La')\to B$ correspond to
$H^0(B,f_*I \ot \TT_B)$. It is easy to check that the map $h_B(\La')\to h^{NC}_B(\La')$ is compatible with the
Kodaira-Spencer map
$$H^0(B,f_*\II \ot\TT_B)\simeq H^0(X_0^{ab},\II\ot f^*\TT_B)\to H^0(X_0^{ab}, \II\ot f^*R^1\rho_*\Eend(\EE)),$$
which is an isomorphism by assumption (b).
It follows that the map $h_B(\La')\to h^{NC}_B(\La')$ is an isomorphism.
\ed



We have the following simple observation.

\begin{prop}\label{fun-factor-prop}
The functor $h^{NC}_B: \NN \to Sets$ factors through $a\NN$.\end{prop}

\Pf. Suppose we have two homomorphims $f_1,f_2:\La'\to \La$ in $\NN$ such that they are conjugate, i.e.,
$f_2=\theta f_1$, where $\th=\th_u$ is an inner automorphism of $\La$: $\th_u(x)=uxu^{-1}$ for some unit $u$ in $\La$.
We have to check that $f_1$ and $f_2$ induce the same map $h(\La')\to h(\La)$. Equivalently, we have to check that
the map $h(\theta):h(\La)\to h(\La)$ is equal to the identity. Note that $\theta_u$ induces an automorphism of the NC-scheme
$X=\Spec(\La)$, which we still denote by $\theta$, and the map $h(\theta)$ sends a right $\OO_{X\times Z}$-module
$E_\La$ to $(\theta\times\id_Z)^*E_\La$. Now we observe that the automorphism $\theta\times\id$ of $X\times Z$
acts trivially on the underlying topological space and is given by the inner
automorphism $\theta_u$ of the structure sheaf $\OO=\OO_{X\times Z}$, associated with $u$ which we view as a global section
of $\OO^*$. Thus, the operation $(\theta\times\id_Z)^*$ is given by
tensoring on the right with the $\OO-\OO$ bimodule $\sideset{_{\theta_u}}{}\OO$ (which is the structure sheaf with the left $\OO$-action twisted by $\theta_u$).

Now we use the general fact that twisting by an inner automorphism does not change an isomorphism class of a bimodule.
Namely, if $M$ is an $R-S$-bimodule and $\theta_u$ is the inner automorphism of $R$ associated with $u\in R^*$, then we have an isomorphism of $(R,S)$-bimodules,
$$M \rTo{\sim} \sideset{_{\theta_u}}{}M: m \mapsto um.$$ 
This construction also works for bimodules over sheaves of rings and an inner automorphism associated with a global unit.
This implies that in our situation the functor $(\theta\times\id_Z)^*$ is isomorphic to identity, 
and our claim follows.
\ed

\begin{rem} In fact, our proof of Proposition \ref{fun-factor-prop} shows a little more. We can enhance $h^{NC}_B$ to
a functor with values in groupoids, by considering the category of the data as in Definition \ref{fun-defi} and isomorphisms between
them. On the other hand, we can consider a $2$-category of algebras in $\NN$ with the usual $1$-morphisms and with
$2$-morphisms between $f_1,f_2:\La'\to\La$ given by $u\in \La^*$ such that $f_2=\theta_u f_1$. Then
the functor $h^{NC}_B$ lifts to a $2$-functor from this $2$-category to the $2$-category of groupoids.
\end{rem}


\begin{thm}\label{non-rep-thm} 
If $\dim B\ge 1$ then for any $d\ge 1$ the functor $h^{NC}_B|_{\NN_d}$ is not representable by an 
NC-nilpotent scheme of degree $d$.
\end{thm}

\Pf. It is enough to consider the case $d=1$. Suppose $h^{NC}_B|_{\NN_1}$ is representable by an NC-nilpotent scheme $X$ of
degree $1$. Then by Proposition \ref{h-smooth-prop}, $X$ is $1$-smooth and $X^{ab}\simeq B$. 
Let $U=\Spec(A)\sub X$ be an affine NC-subscheme corresponding to an open affine subscheme of $B$ of dimension $\ge 1$.
Then $A$ is a $1$-smooth algebra with $\dim A^{ab}\ge 1$, and $h_A$ is a subfunctor of $h^{NC}_B|_{\NN_1}$.
Since the latter functor factors through $a\NN_1$, this would imply that $h_A$ also factors through $a\NN_1$.

It remains to prove that for any $1$-smooth algebra $A$ with $\dim A^{ab}\ge 1$
the functor $h_A$ does not factor through $a\NN_1$. To this end we will give an example of two
conjugate homomorphisms $f_1,f_2:A \to A'$ such that $f_1\neq f_2$.
Set
$$A' =  (A*\C[z,z^{-1}])_{[\![ab]\!]}/\I_3.$$
It is easy to see that $A'$ is $1$-smooth and $(A')^{ab} = A^{ab}\otimes\C[z,z^{-1}]$. 
Therefore, by Lemma \ref{NC-sm-center}(i), the element $z$ is not in the center of $A'$.
Hence, we can take $f_1:A\to A'$ to be the natural homomorphism and set $f_2(a)=zf_1(a)z^{-1}$.
\ed


\section{Representability results}\label{rep-sec}

\subsection{Local representability in $a\NN$}


Kapranov gives the following criterion for a formally smooth functor on $\NN_d$ to be representable by an NC-scheme.

\begin{prop}\label{NC-rep-prop} (\cite[Thm.\ (2.3.5)]{Kapranov})
Let $M$ be a smooth algebraic variety. A formally smooth functor $h: \NN_d \to Sets$ such that $h|_{\Com} = h_M$, is representable by a $d$-smooth NC-scheme if and only if for any pair of central extensions in $\NN_d$, $\La_1\to\La$,
$\La_2\to \La$, the natural map $$h(\Lambda_1\times_\Lambda\Lambda_2) \to h(\Lambda_1)\times_{h(\Lambda)}h(\Lambda_2)$$ is an isomorphism.
\end{prop}


We will prove an analogous representability criterion for affine aNC-schemes.
As in the case of NC-schemes the main idea is to study fibers of the map $h(p):h(\La')\to h(\La)$ for
a central extension 
\begin{equation}\label{central-ext-eq}
0\to I\to \La'\rTo{p} \La\to 0
\end{equation}
(cf. the proof of \cite[Lem.\ (2.3.6)]{Kapranov}).


For $d\ge 1$, let $h:a\NN_d\to Sets$ be a functor such that $h|_{a\NN_{d-1}}$ is representable by $A\in a\NN_{d-1}$.
The key new ingredient we have to use is the following. Given a central extension \eqref{central-ext-eq}
with $\La'\in\NN_d$, $\La\in\NN_{d-1}$, and a homomorphism $f:A\to\La$,
we set
$$U(f):=\{u\in \La^* \ |\ uf(a)u^{-1}=f(a) \forall a\in A\}.$$
Then we have a natural map 
$$\De_f: U(f)\to \Der(A,I)=\Der(A^{ab},I).$$ 
where
\begin{equation}\label{De-f-u-def}
\De_f(u): A\to I: a\mapsto [u,f(a)]_{\La'}u^{-1},
\end{equation}
where for $l_1,l_2\in\La$, we define $[l_1,l_2]_{\La'}\in \La'$ by
\begin{equation}\label{centr-ext-comm-eq}
[l_1,l_2]_{\La'}:=[\wt{l_1},\wt{l_2}],
\end{equation}
where $\wt{l_i}$ is a lifting of $l_i$ to $\La'$. Note that $[u,f(a)]_{\La'}\in I$.

Furthermore, one can check that the image of $\De_f$ depends only on the image of $f$ in 
$\Hom_{a\NN}(A,\La)=h(\La)$. Also, using the fact that $I$ is central we immediately check that
$\De_f$ is a group homomorphism. 
The next result shows that in the case when $h$ itself is representable, the cokernel of $\De_f$ maps bijectively
to $h(p)^{-1}(f)$.

\begin{lem}\label{Der-action-lem} 
Let $A'$ be an NC-nilpotent algebra of degree $d$ such that $A=A'/I_{d+1}A'$. Then for any central
extension \eqref{central-ext-eq}, with $\La'\in\NN_d$ and $\La\in\NN_{d-1}$, and any
algebra homomorphism $f:A'\to \La$ there exists a natural transitive action of the group $\Der(A,I)$ on
the fiber $h_{A'}(p)^{-1}(f)$ of the map $h_{A'}(p):h_{A'}(\La')\to h_{A'}(\La)$, such that the action of $\Der(A,I)$ on
any element of this fiber induces a bijection
$$\coker(\De_f)\rTo{\sim} h_{A'}(p)^{-1}(f).$$
\end{lem} 
 
\Pf . It is well known that the difference between two homomorphisms $A'\to \La'$ lifting $f:A'\to \La$ is a derivation $A'\to I$,
and that this induces a simply transitive action of $\Der(A',I)=\Der(A,I)$ on the set of such liftings. Now assume
that we have two homomorphisms $f'_1,f'_2:A\to \La'$, such that both $p\circ f'_1$ and $p\circ f'_2$ are conjugate to $f$.
Then replacing $f'_1$ and $f'_2$ by conjugate homomorphisms we can assume that $p\circ f'_1=p\circ f'_2=f$.
Now It is easy to see
that if $f'_2$ and $f'_1$ are conjugate by $u\in(\La')^*$ then $u\in U(f)$, and 
the difference $f'_2-f'_1$ is the derivation $a\mapsto [u,f(a)]_{\La'}u^{-1}$.
This establishes the required bijection.
\ed
 

Next, we return to the situation when only $h|_{a\NN_{d-1}}$ is representable.
Recall (see \cite[Prop.\ (1.2.5)]{Kapranov}) 
that for any central extension \eqref{central-ext-eq} there is a natural isomorphism
\begin{equation}\label{La'-square-isom}
\La'\times_\La \La'\rTo{\sim} \La'\times_{\La^{ab}} (\La^{ab}\oplus I): (x,y)\mapsto (x,(x^{ab},y-x)),
\end{equation}
where $x\to x^{ab}$ is the projection $\La'\to \La^{ab}$, and $\La^{ab}\oplus I$ is the trivial commutative
algebra extension of $\La^{ab}$ by $I$ (such that $I^2=0$ and $\La^{ab}$ is a subalgebra).
Let us assume in addition that $h$ commutes with pull-backs by commutative nilpotent extension, so that
$$h(\La'\times_{\La^{ab}} (\La^{ab}\oplus I))\simeq h(\La')\times_{h(\La^{ab})}h(\La^{ab}\oplus I).$$
Combining this with the above isomorphism we get a natural map
\begin{equation}\label{h-main-fiber-map}
h(\La')\times_{h(\La^{ab})}h(\La^{ab}\oplus I)\simeq h(\La'\times_\La \La')\to h(\La')\times_{h(\La)}h(\La').
\end{equation}

Now assume $\La'\in\NN_d$, $\La\in\NN_{d-1}$ and we are given an element $f'\in h(\La')$ lifting $f\in h(\La)$.
Since $h|_{a\NN_{d-1}}\simeq h_A$ we have a natural identification of the fiber of
$h(\La^{ab}\oplus I)\to h(\La^{ab})=\Hom_{alg}(A,\La^{ab})$ over $f^{ab}$ with $\Der(A,I)$.
Thus, for any $D\in \Der(A,I)$ we can consider a pair $(f',f^{ab}+D)$ in the left-hand side of \eqref{h-main-fiber-map}.
Let us define $f'+D\in h(p)^{-1}(f)$, so that $(f',f'+D)$ is the image of $(f',f^{ab}+D)$ under \eqref{h-main-fiber-map}.
In this way we get a map
\begin{equation}\label{add-der-op-eq}
\de_{f'}:\Der(A,I)\to h(p)^{-1}(f): D\mapsto f'+D.
\end{equation}
It is easy to see (by considering $\La'\times_\La \La'\times_\La \La'$) that in this way we get an action of the group $\Der(A,I)$ on
$h(p)^{-1}(f)$.
Note that in the case when $h$ is representable by some $A'\in\NN_d$, this operation is exactly the operation
of adding a derivation $A'\to A\to I$ to a homomorphism $A'\to\La'$. 


Now we can prove the following local aNC version of Proposition \ref{NC-rep-prop}.

\begin{prop}\label{aNC-rep-prop} 
Let $A$ be a $(d-1)$-smooth algebra in $a\NN_{d-1}$, and let $h: a\NN_d \to Sets$, 
be a formally smooth functor such that $h|_{a\NN_{d-1}} \simeq h_A$. Then $h$ is representable by 
a $d$-smooth algebra in $a\NN_d$
if and only if the following two conditions hold.
	\begin{itemize}
		\item[(i)] For any nilpotent extension $\La'\to \La$ with $\La'\in a\NN_d$ and $\La\in\Com$, and any commutative nilpotent extension $\La''\to \La$, the natural map
                $$h(\La'\times_\La \La'')\to h(\La')\times_{h(\La)} h(\La'')$$
                is a bijection.
		\item[(ii)]
For every central extension \eqref{central-ext-eq},  for
any $f'\in h(\La')$ extending $f\in h(\La)$, the map $\de_{f'}$, which is well defined due to condition (i), induces a bijection
$$\coker(\De_f)\rTo{\sim} h(p)^{-1}(f).$$
	\end{itemize}
\end{prop}

\Pf. Assume first that $h$ is representable by $A'\in a\NN_d$. To check condition (i) for $h=h_{A'}$
we first note that since $\La$ and $\La''$
are commutative,
the set $h(\La')\times_{h(\La)} h(\La'')$ can be described as pairs of homomorphisms $f':A\to \La'$ and $f'':A\to \La''$ lifting
the same homomorphism $f:A\to \La$, up to the equivalence replacing $f'$ by a conjugate homomorphism.
Clearly, this is the same as giving a homomorphism $A'\to \La'\times_\La \La''$ up to conjugacy. 
On the other hand, condition (ii) for $h_{A'}$ follows from Lemma \ref{Der-action-lem}.

Now assume that conditions (i) and (ii) hold, and let $A'\to A$ be a $d$-smooth thickening of $A$ 
(it exists by \cite[Prop.\ (1.6.2)]{Kapranov}).
Let $e\in h(A)$ be the family corresponding to the isomorphism $h|_{a\NN_{d-1}} \simeq h_A$.
Since $h$ is formally smooth, there exists an element $e'\in h(A')$ lifting $e$.
Let $h_{A'}\to h$ be the induced morphism of functors.
We already know that it is an isomorphism on $a\NN_{d-1}$, and we claim that it is an isomorphism on $a\NN_d$.
The argument is similar to that of Proposition \ref{h-smooth-prop}. Given $\La'\in \NN_d$, we can fit it into 
a central extension \eqref{central-ext-eq}
with $\La\in \NN_{d-1}$. Then we consider the commutative square
\begin{diagram}
h_{A'}(\La') &\rTo{}& h_{A'}(\La)\\
\dTo{}&&\dTo{}\\
h(\La') &\rTo{} & h(\La)
\end{diagram}
Since $h_{A'}(\La)\simeq h_A(\La)\simeq h(\La)$, we know that the right vertical arrow is an isomorphism.
Also, both horizontal arrows are surjective. Let us fix a homomorphism $f\in h_A(\La)$, and its lifting $f'\in h_{A'}(\La')$.
As we have seen in Lemma \ref{Der-action-lem}, the fiber of the top horizontal arrow over $f$ is identified with
$\coker(\De_f)$. The same is true for the fiber of the bottom horizontal arrow over $f$, by condition (ii).
It remains to observe that both isomorphism are induced by the operation \eqref{add-der-op-eq}
of adding a derivation in $\Der(A,I)$, which is compatible with morphisms of functors on $a\NN_d$, extending $h_A$ on
$a\NN_{d-1}$. Thus, the left vertical arrow induces an isomorphism between the fibers of the horizontal arrows over $f$.
Since $f$ was arbitrary, we deduce that the left vertical arrow is an isomorphism.
\ed


\begin{rem} All the fiber products of algebras above are taken in $\NN_d$. 
Fiber products in $a\NN_d$ usually do not exist (unless one of the factors is commutative).\end{rem}

\subsection{Local representability of the functor of NC-families by an aNC scheme}

Assume we are in the situation of Sec.\ \ref{versal-sec}. By Proposition \ref{fun-factor-prop}, 
we can view $h^{NC}_B$ as a functor
on the category $a\NN$,
Our main goal is to prove the local representability of the corresponding functor $h^{NC}_B|_{a\NN_d}$ by a $d$-smooth
NC-algebra.

\begin{thm}\label{loc-rep-thm} 
Assume that the base $B$ of an excellent family is affine. Then for every $d\ge 0$ the functor
$h^{NC}_B|_{a\NN_d}$ is representable by a $d$-smooth thickening of $B$. Hence the functor $h_B^{NC}$ is representable
by a NC-smooth thickening of $B$.
\end{thm}







The proof will proceed by induction on $d$. 
We need two technical Lemmas (the second of which is a noncommutative 
extension of Lemma \ref{comm-End-lem}).

\begin{lem}\label{com-KS-lem}
Assume that $h^{NC}_B|_{a\NN_{d-1}}$ is representable by $A\in\NN_{d-1}$. Then for any central extension
\eqref{central-ext-eq} with $\La\in a\NN_{d-1}$, $\La'\in a\NN_d$, and any homomorphism $f:A\to \La$, there is a commutative
square
\begin{equation}\label{KS-comm-sq}
\begin{diagram}
U(f)&\rTo{\De_f}& \Der(A^{ab},I)\\
\dTo{}&&\dTo{-KS}\\
\Aut(E_\La)&\rTo{\de_0}& H^1(\Spec(\La^{ab})\times Z,\End(E^{ab})\otimes I)
\end{diagram}
\end{equation}
Here $\De_f$ is given by \eqref{De-f-u-def}; 
$E_\La=E_f$ is the family in $h^{NC}_B(\La)$ induced by $f$; the map $KS$
is induced by the Kodaira-Spencer map; and the homomorphism $U(f)\to \Aut(E_f)$ associates with $u\in\La^*$ 
an automorphism of $E_f$ induced by the left multiplication by $u$ on $\La$.
The map $\de_0$ is the connecting map associated with the exact sequence of sheaves
\eqref{Aut-ex-seq},
where $E_{\La'}$ is a vector bundle over $\Spec(\La')\times Z$ lifting $E_\La$.
In particular, in this situation $\de_0$ is a group homomorphism.
\end{lem}

\Pf . We are going to compute the maps in the square \eqref{KS-comm-sq} using local trivializations.
Let us denote by $\EE^{ab}$ the original family over $B\times Z$, and let $\EE$ be the family over $\Spec(A)\times Z$ corresponding to the element $\id_A\in h_A(A)\simeq h^{NC}_B(A)$.
We denote by $f^{ab}$ the homomorphism $A^{ab}\to\La^{ab}$ induced by $f$ and the corresponding morphism
of affine schemes $\Spec(\La^{ab})\to \Spec(A^{ab})=B$. Note that by Proposition \ref{h-smooth-prop}, we have
an isomorphism $E^{ab}=(f^{ab}\times \id)^*\EE^{ab}$.

\medskip

\noindent
{\bf Step 1}. Computation of $\de_0: \Aut(E_f)\to H^1(\Spec(\La^{ab})\times Z,\Eend(E^{ab})\otimes p_1^*\II)$. 

Let us fix an open affine covering $(U_i)$ of $\Spec(\La^{ab})\times Z$ such that $E_{f'}$ is trivial over $U_i$. 
Then, given an automorphism $\a\in\Aut(E_f)$, over $U_i$ we can lift $\a$ to an automorphism $\a_i$ of $E_{\La'}$. 
Now over $U_i\cap U_j$ the endomorphism $\a_i^{-1}\a_j-\id$
of $E_{\La'}$ factors through the kernel
of the projection $E_{f'}\to E_f$, i.e., $E^{ab}\ot p_1^*\II$. This gives the Cech $1$-cocycle with values in 
$\Eend(E^{ab})\ot p_1^*\II$, representing the class $\de_0(\a)$.

\medskip

\noindent
{\bf Step 2}. Computation of the KS-map 
\begin{equation}\label{KS-map}
\Der(A^{ab},I)\to H^1(\Spec(\La^{ab})\times Z,\Eend(E^{ab})\otimes p_1^*\II).
\end{equation}

Note that we have an identification
$$\Der(A^{ab},I)\simeq H^0(B,\TT_B\ot f^{ab}_*\II).$$ 
Let us fix trivializations $\varphi^{ab}_i:\OO^n\to \EE^{ab}$ over an affine open covering $(U_i)$ of $B\times Z$, 
and let $g^{ab}_{ij}=(\varphi^{ab}_i)^{-1}\varphi^{ab}_j\in \Mat_n(\OO(U_i\cap U_j))$ be the corresponding transition functions. Then
to a vector field $v$ on $B$ with values in $f^{ab}_*\II$ 
the KS-map associates the Cech $1$-cocycle $\varphi^{ab}_i v(g^{ab}_{ij})(g^{ab}_{ij})^{-1} (\varphi^{ab}_i)^{-1}$ on $B\times Z$
with values in $\Eend(\EE^{ab})\ot p_1^*f^{ab}_*\II$.

We also need to calculate the image of this class under the isomorphism induced by the projection formula
\begin{align*}
&H^1(B\times Z,\Eend(\EE^{ab})\otimes p_1^*f_*\II)\rTo{\sim}
H^1(B\times Z, (f\times\id)_*((f\times \id)^*\Eend(\EE^{ab})\otimes p_1^*\II))\simeq\\
&H^1(\Spec(\La^{ab})\times Z,\Eend(E^{ab})\otimes p_1^*\II).
\end{align*}
To this end we note that the morphism $f^{ab}\times \id:\Spec(\La^{ab})\times Z\to B\times Z$ is affine, and
so $\wt{U}_i:=(f^{ab}\times\id)^{-1}(U_i)$ is an affine open covering of $\Spec(\La^{ab})\times Z$, over which we have
the induced trivializations of $E^{ab}=(f^{ab}\times\id)^*\EE^{ab}$, which we still denote by $\varphi^{ab}_i$.
Now it is easy to see that the corresponding Cech 1-cocycle on $\Spec(\La^{ab})\times Z$
with values in $\Eend(E^{ab})\otimes I$ is given by
$$\varphi^{ab}_i v(g^{ab}_{ij}) f^{ab}(g^{ab}_{ij})^{-1} (\varphi^{ab}_i)^{-1},$$ 
where we denote still by $f^{ab}:\OO(U_i\cap U_j)\to \OO(\wt{U}_i\cap \wt{U}_i)$ the homomorphism induced by $f^{ab}$,
and also extend $v$ to a derivation $\OO(U_i\cap U_j)\to p_1^*\II(\wt{U}_i\cap\wt{U}_j)$.

\medskip

\noindent
{\bf Step 3}. Now we can check the commutativity of the square \eqref{KS-comm-sq}

We start by choosing an affine open covering $(U_i)$ of $B\times Z$ and trivializations of $\EE^{ab}$ over $U_i$.
Then we can lift these trivializations to some trivializations $\varphi_i:\OO^n_{\Spec(A)\times Z}|_{U_i}\to \EE$ 
(see Lemma \ref{triv-bun-lem}).
We denote by $g_{ij}$ the corresponding transition functions in $\GL_n(\OO_{\Spec(A)\times Z}(U_i\cap U_j))$.

By definition, $\De_f(u)$ is the derivation 
$$v(a) = [u, f(a)]_{\La'} u^{-1}=[\tilde u,\wt{f(a)}]\tilde u^{-1},$$
where $\tilde u, \wt{f(a)} \in \La'$ are some lifts of $u$ and $f(a)$ 
(note that $\Der(A,I)=\Der(A^{ab},I)$).
Hence, $KS(\Delta_f(u))$ is represented by the 1-cocycle 
\begin{equation}\label{KS-Delta-eq}
\varphi_i [\tilde u, \wt{f(g_{ij})}]_{\La'}\tilde u^{-1}\wt{f(g_{ij})}^{-1}\varphi_i^{-1} = \varphi_i(\tilde u \wt{f(g_{ij})}
\tilde u^{-1}\wt{f(g_{ij})}^{-1} - \id)\varphi_i^{-1}.
\end{equation}

As in Step 2, we have the induced affine open covering $\wt{U}_i$ of $\Spec(\La^{ab})\times Z$, and the induced trivializations
$\psi_i$ of $E_f$ over $\wt{U}_i$. Let us choose a lifting $E_{\La'}$ of $E_f$ to a vector bundle over $\Spec(\La')\times Z$
(it exists by formal smoothness of $h^{NC}_B$),
and liftings $\psi'_i$ of $\psi_i$ to trivializations of $E_{\La'}$ over $\wt{U}_i$ (see Lemma \ref{triv-bun-lem}).
Note that we have $\psi_i^{-1}\psi_j = f(g_{ij})$, and hence $(\psi'_i)^{-1}\psi'_j$ provide liftings $\wt{f(g_{ij})}\in\La'$ of $f(g_{ij})$. 
The image of $u \in U(f)$ in $\Aut(E_f)$ can be represented over $\wt{U}_i$ as $\psi_i u \psi_i^{-1}$, where we view $u$
as the corresponding operator
of the left multiplication by $u$ (note that these operators are compatible on intersections because
$u\cdot f(g_{ij}) = f(g_{ij})\cdot u$, due to the inclusion $u \in U(f)$). 
Using the lifting $\tilde u \in \La'$ of $u$
we get local automorphisms of $E_{f'}$ over $\wt{U}_i$,
$\alpha_i = \psi'_i \tilde u(\psi'_i)^{-1}$. 
Then 
$$\delta_0(\alpha) = \alpha_i^{-1}\alpha_j - \id = (\psi'_i\tilde u^{-1} (\psi'_i)^{-1})(\psi'_j\tilde u \tilde \psi_j^{-1}) - \id = 
\psi'_i(\tilde u^{-1} \wt{f(g_{ij})}\tilde u \wt{f(g_{ij})}^{-1} -\id)(\psi'_i)^{-1}.$$ 
Comparing this with \eqref{KS-Delta-eq} we see that
$$\de_0(\a)=KS(\Delta_f(u^{-1}))=KS(-\De_f(u))=-KS(\De_f(u)).$$
\ed

\begin{lem}\label{main-Aut-lem} 
Assume that $h^{NC}_B|_{a\NN_d}$ is representable by $A\in a\NN_d$, so $h^{NC}_B|_{a\NN_d}\simeq h_A$. 
Then for every $d$-nilpotent algebra $\La$ and every homomorphism $f:A\to \La$, the induced homomorphism 
$U(f)\to \Aut(E_f)$ is an isomorphism. Here $E_f$ represents the family in $h^{NC}_B(\La)$ induced by $f$. 
\end{lem}

\Pf . 
We will prove the assertion by induction on $d'\le d$ such that $\La$ is $d'$-nilpotent. For $d'=0$, i.e., when $\La$ is commutative,
we have $U(f)=\La^*$ and the assertion follows from Lemma \ref{comm-End-lem}.

Next, we have to see that both groups fit into the same exact sequences, when $\La'$ is a central extension of $\La$ by $I$.
Namely, if $f':A\to \La'$ is a homomorphism lifting $f$, then by Lemma \ref{com-KS-lem},
we have a morphism of exact sequences
\begin{diagram}
1 &\rTo{}& 1+I&\rTo{} &U(f')&\rTo{}& U(f)&\rTo{\De_f}& \Der(A^{ab},I)\\
&&\dTo{\id}&&\dTo{}&&\dTo{}&&\dTo{-KS}\\
1 &\rTo{}& 1+I&\rTo{} &\Aut(E_{f'})&\rTo{}& \Aut(E_f)&\rTo{\de_0}& H^1(\Spec(\La^{ab})\times Z,\Eend(E^{ab})\otimes p_1^*\II)
\end{diagram}
Note that the map $KS$ is an isomorphism. Since the map $U(f)\to\Aut(E_f)$ is an isomorphism by the
induction assumption, we deduce that $U(f')\to \Aut(E_{f'})$ is also an isomorphism.
\ed

\noindent
{\it Proof of Theorem \ref{loc-rep-thm}}.
By Proposition \ref{h-smooth-prop}, we know that the assertion is true for $d=0$. 
Now, assuming that the functor $h^{NC}_B|_{a\NN_{d-1}}$ is representable, we will apply Proposition \ref{aNC-rep-prop} to prove
that $h^{NC}_B|_{a\NN_d}$ is representable. It suffices to check conditions (i) and (ii) of this Proposition.
To prove condition (i) assume that $\La'\to \La$ and $\La''\to \La$ and nilpotent extensions with $\La,\La''\in\Com$.
To see that the map
$$h(\La'\times_\La \La'')\to h(\La')\times_{h(\La)} h(\La'')$$
is a bijection, we construct (as in \cite[Lem.\ (5.4.4)]{Kapranov}) the inverse map as follows.
Starting with families $\EE_{\La'}$ and $\EE_{\La''}$ over $\La'$ and $\La''$, and choosing an arbitrary
isomorphism of the induced families over $\La$, we define the family over $\La'\times_\La \La''$ as
the fibered product $\EE_{\La'}\times_{\EE_{\La}}\EE_{\La''}$. One has to check that the result does not depend on a choice
of isomorphism of families over $\La$ (this may fail in general, but works for commutative $\La''$).
Note that different choices differ by an automorphism of $\EE_{\La}$, so it is enough to see that any such automorphism can be
lifted to an automorphism of $\EE_{\La''}$. But this follows immediately from Lemma \ref{comm-End-lem}.

Next, let us check condition (ii). Given a central extension \eqref{central-ext-eq} with $\La'\in \NN_d$, $\La\in \NN_{d-1}$,
and a family $(f^{ab},E_\La,\phi)$ in $h^{NC}_B(\La)$,  
then choosing a lifting $E_{\La'}$ to a family over $\La'$,
from the corresponding exact sequence of sheaves of groups \eqref{Aut-ex-seq} we get a connecting map
$$\de_0:\Aut(E_\La)\to H^1(X^{ab}\times Z,\Eend(E_{\La^{ab}})\otimes p_1^*\II).$$
Furthermore, by Lemma \ref{com-KS-lem}, $\de_0$ is actually a group homomorphism (and the source of this map acts trivially on
the target).
Thus, from the formalism of nonabelian cohomology applied to the abelian extension of sheaves of groups \eqref{GL-ex-seq}
we get that different liftings of $E_\La$ to a family over $\La'$ form a
principal homogeneous space over $\coker(\de_0)$ (see Sec.\ref{nonab-H1-sec}).
Note that by Lemma \eqref{main-Aut-lem}, we have an isomorphism $U(f)\simeq \Aut(E_\La)$, where
$f:A\to \La$ is the homomorphism giving $E_\La$. 
Thus, by Lemma \ref{com-KS-lem}, we can identify $\coker(\de_0)$ with $\coker(\De_f)$.
Thus, to prove condition (ii), it remains to check that the two actions of $\Der(A,I)$
on the set of liftings of $E_\La$ are the same (the one coming from the formalism of non-abelian cohomology,
and the other one given by the map \eqref{add-der-op-eq}).

To this end we use the computation of the Kodaira-Spencer map \eqref{KS-map} using local trivializations.
Namely, we choose trivializations of the universal bundle $\EE$ over an open covering of $\Spec(A)\times Z$,
and denote by $g_{ij}$ the corresponding transition functions, so that $f(g_{ij})$ are the transition functions for $E_\La$.
Then, in the notation of Lemma \ref{com-KS-lem}, a derivation $v\in\Der(A,I)=\Der(A^{ab},I)$ 
gives rise to the Cech $1$-cocycle
$$\varphi_i v(g_{ij})f(g_{ij})^{-1}\varphi_i^{-1}$$
on $\Spec(\La^{ab})\times Z$ with values in $\Eend(E^{ab})\ot p_1^*\II$.
The corresponding $f(g_{ij})$-twisted $1$-cocycle with values in $\Mat_r(\OO)\ot p_1^*\II$
is $(v(g_{ij})f(g_{ij})^{-1})$. Now by definition, the action of $v$ on the set of liftings of $f(g_{ij})$ to
a $1$-cocycle with values in $\GL_r(\OO_{\Spec(\La')\times Z})$ sends $(\wt{g}_{ij})$ to
\begin{equation}\label{Der-action-v-g-cocycle-eq}
((1+v(g_{ij})f(g_{ij})^{-1})\cdot \wt{g}_{ij}=(\wt{g}_{ij}+v(g_{ij})).
\end{equation}
On the other hand, from $v$ we get a homomorphism $f^{ab}+v:A\to \La^{ab}\oplus I$, and hence,
the $1$-cocycle $(f^{ab}+v)(g_{ij})$ with values in $\GL_r(\OO_{\Spec(\La^{ab}\oplus I)\times Z})$ lifting $f^{ab}(g_{ij})$.
Hence, a lifting $\wt{g}_{ij}$ of $f(g_{ij})$ together with $v$ defines a $1$-cocycle
$$(\wt{g}_{ij}, (f^{ab}+v)(g_{ij}))$$
with values in $\GL_r(\OO_{\Spec(\La'\times_{\La^{ab}}(\La^{ab}\oplus I))\times Z})$.
It remains to observe that under the isomorphism \eqref{La'-square-isom} it corresponds to the $1$-cocycle
$$(\wt{g}_{ij}, \wt{g}_{ij}+v(g_{ij}))$$
with values in $\GL_r(\OO_{\Spec(\La'\times_{\La}\La')\times Z})$, which has \eqref{Der-action-v-g-cocycle-eq}
as the same second component.
\ed


\subsection{Nonabelian hypercohomology}\label{nonab-hypercoh}

We will use below the following simple generalization of nonabelian $H^1$.
Let $\GG$ be a sheaf of groups over a topological space $X$, and let $\EE$ be a sheaf of sets,
equipped with a $\GG$-action. We view a pair $\GG\car \EE$ as a generalization of a length $2$ complex.

For an open covering $\UU=(U_i)_{i\in I}$ of $X$, we define the set
of $1$-cocycles over $\UU$ for the pair $\GG\car \EE$:
\begin{align}
Z^1(\UU,\GG\car\EE):=&\{(g_{ij}\in\GG(U_{ij}))_{i,j\in I}, (e_i\in\EE(U_i))_{i\in I} \ |\ g_{ii}=1, \ g_{ij}g_{ji}=1, \ \nonumber\\
& g_{ij}g_{jk}=g_{ik}, \  e_i=g_{ij}(e_j) \},
\end{align}
where as usual we denote $U_{ij}=U_i\cap U_j$, $U_{ijk}=U_i\cap U_j\cap U_k$ (and the restrictions to appropriate
intersections are assumed).
Two $1$-cocycles over $\UU$, $(g_{ij},e_i)$ and $(\wt{g}_{ij},\wt{e}_i)$ are called cohomologous if for some
collection $h_i\in \GG(U_i)$ we have
$$\wt{g}_{ij}=h_ig_{ij}h_j^{-1}, \ \ \wt{e}_i=h_i(e_i).$$
It is easy to see that this defines an equivalence relation on $Z^1(\UU,\GG\car\EE)$, and we denote by
$\hH^1(\UU,\GG\car\EE)$ the corresponding set of equivalence classes. Passing to the limit over all open coverings $\UU$,
we get the {\it nonabelian hypercohomology set} $\hH^1(X,\GG\car\EE)$.

These sets are natural: if we have a homomorphism of sheaves of groups $\GG_1\to \GG_2$ and the
compatible map of sheaves of set $\EE_1\to \EE_2$, then we get the induced map 
$$\hH^1(X,\GG_1\car\EE_1)\to \hH^1(X,\GG_2\car\EE_2).$$
Also, sending $(g_{ij},e_i)$ to $g_{ij}$ defines a projection to the usual nonabelian $H^1$,
$$\hH^1(X,\GG\car\EE)\to H^1(X,\GG).$$

Recall that $H^1(X,\GG)$ classifies isomorphism classes of $\GG$-torsors. Similarly, the set $\hH^1(X,\GG\car\EE)$
can be identified with the isomorphism classes of pairs $(P,e)$, where $P$ is a $\GG$-torsor, and $e$ is a global
section of the twisted sheaf $\EE_P=P\times_{\GG} \EE$. 

Next, we have the following analog of the connecting homomorphism $H^1\to H^2$.
Assume that we have an abelian extension of sheaves of groups
$$1\to \AA_0\to \GG'\rTo{p} \GG\to 1$$
over $X$, and sheaves of sets $\EE'$ and $\EE$, where $\GG'$ (resp., $\GG$) acts on $\EE'$ (resp., $\EE$).
Further, assume that we have a sheaf of abelian groups $\AA_1$ acting freely on $\EE'$, and an identifcation
$\EE=\EE'/\AA_1$. We denote this action as $a_1+e'$, where $a_1\in\AA_1$, $e'\in \EE'$. 
We require the following compatibilities between these data.
First, the projections $p:\EE'\to \EE$ and $p:\GG'\to \GG$ should be compatible with the actions (of $\GG'$ on $\EE'$
and of $\GG$ on $\EE$). Note that this implies that there is an action of $\GG'$ on $\AA_1$, compatible
with the group structure on $\AA_1$, such that
$$g'(a_1+e')=g'(a_1)+g'(e').$$
Secondly, we require that the subgroup $\AA_0\sub \GG'$ acts trivially on $\AA_1$, so that there is an induced action
of $\GG$ on $\AA_1$, such that the above formula becomes
$$g'(a_1+e')=p(g)(a_1)+g'(e').$$
In particular, for $g'=a_0\in\AA_0$, we get
\begin{equation}\label{a0-a1-e'-eq}
a_0(a_1+e')=a_1+a_0(e').
\end{equation}
For $e'\in \EE'$ and $a_0\in\AA_0$, let us define $d_{e'}(a_0)\in\AA_1$ from the equation
$$a_0(e')=d_{e'}(a_0)+e'$$
(this is possible since $a_0$ acts trivially on $\EE$).
Furthermore, \eqref{a0-a1-e'-eq} easily implies that $d_{a_1+e'}(a_0)=d_{e'}(a_0)$, so we have
a well defined map of sheaves
$$\EE\times \AA_0\to \AA_1: (e,a_0)\mapsto d_e(a_0),$$
compatible with the group structures in $\AA_0$ and $\AA_1$, such that
$$a_0(e')=d_{p(e')}(a_0)+e'.$$
In particular, for every section $e$ of $\EE$ over an open subset $U\sub X$ 
we have a complex of abelian groups over $U$, $(\AA_\bullet,d_e)$.
Note that $\GG$ acts on $\AA_0$ (via adjoint action $\Ad(g)$), $\AA_1$ and $\EE$, and we have
\begin{equation}\label{g-d-e-a0-eq}
g(d_e(a_0))=d_{g(e)}(\Ad(g)a_0).
\end{equation}

Now assume we have a class $c\in \hH^1(X,\GG\car\EE)$ represented by
a Cech $1$-cocycle $(g_{ij},e_i)$. Let $g=(g_{ij})$ be the induced class
in $H^1(X,\GG)$. We have the corresponding twisted
sheaves $\AA_0^g$, $\AA_1^g$, and \eqref{g-d-e-a0-eq} implies that the $d_{e_i}$'s glue into a global differential
$$d_e:\AA_0^g\to \AA_1^g.$$ 
We are going to define an obstruction class $\de_1(c)$ with values in
$$\hH^2(X,(\AA_\bullet^g,d_e)),$$
such that it vanishes if and only if $(g_{ij},e_i)$ can be lifted to a class in $\hH^1(X,\GG'\car\EE')$.
Namely, by making the covering small enough, we can assume that 
$$g_{ij}=p(g'_{ij}), \ g'_{ij}\in \GG'(U_{ij}), \ \ e_i=p(e'_i), \ e'_i\in \EE'(U_i).$$
Then we have well defined elements $a_{0,ijk}\in\AA_0(U_{ijk})$ and $a_{1,ij}\in \AA_1(U_{ij})$,
such that
$$g'_{ij}g'_{jk}=a_{0,ijk}g'_{ik},$$
$$g'_{ij}(e'_j)=a_{1,ij}+e'_i.$$
It is easy to check that $(a_{0,ijk}, a_{1,ij})$ satisfy the equations
$$a_{0,ijk}+a_{0,ikl}=\Ad(g_{ij})a_{0,jkl}+a_{0,ijl}, \ \ a_{1,ij}+g_{ij}(a_{1,jk})=d_{e_i}(a_{0,ijk})+a_{1,ik},$$
which exactly means that we get a $2$-cocycle $\de_1(g_{ij},e_i)$ with values in $(\AA_\bullet^g,d_e)$.

One can check that this construction gives a well defined element  
$\de_1(c)\in\hH^2(X,(\AA_\bullet^g,d_e))$. Namely, a different choice of liftings $g'_{ij}\mapsto a_{0,ij}g'_{ij}$,
$e'_i\mapsto a_{1,i}+e'_i$ would lead to adding the coboundary of $(a_{0,ij}, a_{1,i})$ to 
the twisted $2$-cocycle $(a_{0,ijk},a_{1,ij})$. On the other hand, changing $(g_{ij},e_i)$ to
$(h_ig_{ij}h_j^{-1},h_i(e_i))$ would lead to a different presentation of the twisted sheaves $\AA_\bullet^g$,
so that the action of $h_i$ glues into isomorphism between two presentations. Our $2$-cocycles $\de_1(g_{ij},e_i)$
and $\de_1(h_ig_{ij}h_j^{-1},h_i(e_i))$ correspond to each other under this isomorphism.

Next, let us assume that a class $c\in \hH^1(X,\GG\car\EE)$ is lifted to a class $c'\in \hH^1(X,\GG'\car\EE')$.
(More precisely, we need to fix the corresponding pair $(P',e')$ where $P'$ is $\GG'$-torsor and $e'$ is a
global section of $\EE'_{P'}$.)
Let $g\in H^1(X,\GG)$ be the image of $c$. We define the following subgroup in $H^0(X,\GG^g)$:
$$\hH^0(X,\GG,c):=\{ (\a_i\in \GG(U_i)) \ |\ \a_i=g_{ij}\a_jg_{ij}^{-1}, \ \a_i(e_i)=e_i\},$$
where $(g_{ij},e_i)$ is a Cech representative of $c$. 
We have a natural connecting map (depending on a choice of $c'$)
$$\de_0:\hH^0(X,\GG,c)\to \hH^1(X,(\AA_\bullet^g,d_e)),$$
defined as follows. We can assume $(g_{ij},e_i)$ comes from a Cech representative $(g'_{ij},e'_i)$ for $c'$.
Let $\a=(\a_i)$ be an element in $\hH^0(X,\GG,c)$. We can assume
that each $\a_i$ can be lifted to $\a'_i\in\GG'(U_i)$. Then we have
$$\a'_i\cdot a_{0,ij}=g'_{ij}\a'_j(g'_{ij})^{-1}, \ \ \a'_i(a_{1,i}+e'_i)=e'_i,$$
for uniquely defined $a_{0,ij}\in \AA_0(U_{ij})$, $a_{1,i}\in \AA_0(U_i)$.
It is easy to check that the following equations are satisfied:
\begin{equation}\label{g-twisted-1cocycle-a0-a1-eq}
a_{0,ij}+\Ad(g_{ij})(a_{0,jk})=a_{0,ik}, \ \ d_{e_i}(a_{0,ij})=a_{1,i}-g_{ij}(a_{1,j}),
\end{equation}
which mean that $(a_{0,ij},a_{1,i})$ define a $1$-cocycle with values in $(\AA_\bullet^g,d_e)$.
We set $\de_0(\a_i)$ to be the class of this $1$-cocycle.
As in Sec.\ \ref{nonab-H1-sec}, one can check that $\a\mapsto \de_0(\a^{-1})$ is a crossed homomorphism,
i.e., equation \eqref{de0-crossed-eq} is satisfied.

Next, we have a natural surjective map (depending on $c'$)
\begin{equation}\label{H1-Lc-map-eq}
\hH^1(X,(\AA_\bullet^g,d_e))\to L_c,
\end{equation}
where $L_c\sub \hH^1(X,\GG'\car\EE')$ is the set of liftings of $c$.
Namely, given a twisted Cech $1$-cocycle with values in $(\AA_\bullet^g,d_e)$, $(a_{0,ij},a_{1,i})$, so that
equations \eqref{g-twisted-1cocycle-a0-a1-eq} are satisfied,
and a representative $(g'_{ij},e'_i)$ of $c'$ we get a new lifting $(a_{0,ij}g'_{ij},a_{1,i}+e'_i)$.
Furthermore, as in Sec.\ \ref{nonab-H1-sec}, 
we can identify the fibers of \eqref{H1-Lc-map-eq} with the orbits of the twisted action of
$\hH^0(X,\GG,c)$ on $\hH^1(X,(\AA_\bullet^g,d_e))$, which is defined similarly to \eqref{twisted-H0-action-eq}.
In particular, in the case when the usual action of $\hH^0(X,\GG,c)$ on $\hH^1(X,(\AA_\bullet^g,d_e))$ is trivial
(or equivalently, $\de_0$ is a group homomorphism), these orbits are simply the cosets for the image of $\de_0$.

\subsection{Families of representations of quivers}

Now we are going to consider families of representations of quivers (without relations).
Let $Q$ be a finite quiver with the set of vertices $Q_0$ and the set of arrows $Q_1$.
We denote by $h,t:Q_1\to Q_0$ the maps associating with an arrow its head and tail.

As in \cite{Toda}, we can consider representations of $Q$ over an NC-scheme $X$. 
Such a representation is a collection of vector bundles $(\VV_v)_{v\in Q_0}$ over $X$, and a collection
of morphisms $e_a:\VV_t(a)\to \VV_{h(a)}$, for each $a\in Q_1$. 

With a collection $\VV=(\VV_v)_{v\in Q_0}$ 
of vector bundles over $X$ we associate a triple of sheaves of groups on the underlying
topological space of $X$,
$$\GG(\VV):=\prod_v \und{\Aut}(\VV_v), \ \ \EE_0(\VV):=\prod_v \Eend(\VV_v), \ \ 
 \EE_1(\VV):=\prod_a \Hhom(\VV_{t(a)},\VV_{h(a)}).$$
Note that there is a natural action of $\GG(\VV)$ on $\EE_1(\VV)$ given by 
$$(g_v)\cdot (\phi_a)=(g_{h(a)}\phi_a g_{t(a)}^{-1}).$$
In the case of trivial bundles $\VV_v=\OO^{n_v}$, for a dimension vector $n_\bullet$, 
we denote these sheaves as $\GG(n_\bullet)$, $\EE_0(n_\bullet)$ and $\EE_1(n_\bullet)$. When we want
to stress the dependence on the NC-scheme $X$ we write $\GG(n_\bullet,X)$, etc.

A structure of a representation of $Q$ on $\VV$ is given by a global section $e=(e_a)$ of $\EE_1(\VV)$.
For such a structure $e$ we can build a 2-term complex 
$$\EE_\bullet(\VV,e): \EE_0(\VV)\rTo{d_f}\EE_1(\VV),$$ 
where the differential is given by $d_e(\phi_v)=\phi_{h(a)}e_a-e_a\phi_{t(a)}$.
Note that $\HH^0\EE_\bullet(\VV,e)$ is precisely the sheaf of endomorphisms of $(\VV,e)$ as a representation of $Q$.


Let $(\VV,e)$ be a representation of $Q$ over $X$.
Over some open affine covering $\UU=(U_i)$ of $X$ we can choose a trivialization 
$\varphi_i=(\varphi_{v,i}):\bigoplus_v\OO^{n_v}_{U_i}\to \bigoplus_v\VV_v|_{U_i}$.
Then over each $U_i$ we have morphisms 
$$e_{a,i}:=\varphi^{-1}_{h(a),i} e_a\varphi_{t(a),i}\in \Mat_{n_{t(a)}\times n_{h(a)}}(\OO(U_i))=\EE_1(n_\bullet)(U_i),$$ 
and over intersections $U_i\cap U_j$ we have transition functions 
$$g_{ij}=(g_{v,ij})=\varphi_i^{-1}\varphi_j\in\prod_v \GL_{n_v}(\OO(U_i\cap U_j))= \GG(n_\bullet)(U_i\cap U_j).$$
One immediately checks that $(g_{ij},e_{a,i})$ defines a Cech $1$-cocycle 
with values in the pair $\GG(n_\bullet)\car \EE_1(n_\bullet)$ (see Sec.\ \ref{nonab-hypercoh}). 
Furthermore, a different choice of trivializations $(\varphi_i)$ leads to a cohomologous cocycle,
so we have a well defined element of $\hH^1(X,\GG(n_\bullet)\car \EE_1(n_\bullet))$.
One can easily check that in this way we get a bijection between the latter nonabelian hypercohomology set
and the set of isomorphism classes of representations $(\VV,e)$ of $Q$, such that the underlying vector bundle has
dimension vector $n_\bullet$.

For a central extension \eqref{central-ext-eq} we have an abelian extension of sheaves of groups 
\begin{equation}\label{GL-ex-seq-bis}
1\to \EE_0(n_\bullet,\OO_{X^{ab}})\ot\II\to \GG(n_\bullet,\OO_{X'})\to \GG(n_\bullet,\OO_{X})\to 1
\end{equation}
where $X=\Spec(\La)$, $X'=\Spec(\La')$, $\II\sub\OO_{X'}$ is the ideal sheaf associated with $I$,
and an exact sequence of abelian groups
$$0\to \EE_1(n_\bullet)\ot\II\to \EE_1(n_\bullet,X')\to \EE_1(n_\bullet,X)\to 0,$$
compatible with the actions of the groups from \eqref{GL-ex-seq-bis}. 
From Sec.\ \ref{nonab-hypercoh} we get that the obstacle to lifting a 
representation $(\VV,e)$ of $Q$ over $\Spec(\La)$ to
a representation of $Q$ over $\Spec(\La')$ is an element of the hypercohomology
$\hH^2(X^{ab}, \EE_\bullet(\VV,e)\ot\II)$. But the latter group $\hH^2$ fits into the exact sequence
$$\ldots\to H^1(X^{ab}, \EE_1(\VV)\ot \II)\to \hH^2\to
H^2(X^{ab}, \EE_0(\VV)\ot\II)\to \ldots$$
Since $X^{ab}$ is an affine scheme, we deduce that our $\hH^2$ vanishes. Thus, the functor of families
of $Q$-representations on $\NN$ is formally smooth.

\begin{defi}\label{quiver-KS}
With a representation $(\VV,e)$ of $Q$ over a commutative scheme $B$ we associate the KS-map, which is a 
morphism of coherent sheaves on $B$,
\begin{equation}\label{quiver-KS-map}
KS:\TT_B\to \HH^1\EE_\bullet(\VV,e),
\end{equation}
defined as follows. Locally we can choose trivializations $\varphi:\bigoplus_v \OO^{n_v}\to \bigoplus_v \VV_v$ and
set for a local derivation $v$ of $\OO_B$,
$$KS(v):=\varphi v(\varphi^{-1} e_a\varphi) \varphi^{-1} \mod \im(d_e)\in \EE_1(\VV,f))/\im(d_e).$$
It is easy to check that a change of a local trivialization leads to an addition of a term in $\im(d_e)$, so the map $KS$
is well defined.
\end{defi}

This definition is motivated by the fact that in the case when
$B=\Spec(k)$ is the point and $(V,e)$ is a $Q$-representation over $k$,
the space $H^1\EE_\bullet(V,e)$ is isomorphic to $\Ext^1((V,e),(V,e))$ (see \cite[Cor.\ 1.4.2]{Brion}), which is
the tangent space to deformations of $(V,e)$ as a $Q$-representation.

Now let us fix a family $(\VV^{ab},e^{ab})$ of representations of $Q$ over a smooth commutative base scheme $B$.
We have the following analog of Definition \ref{versal-fam-defi}.

\begin{defi}\label{versal-quiver-def}
We say that $(\VV^{ab},e^{ab})$ is an {\it excellent family} of representations of $Q$ if
\begin{enumerate}
	\item[(a)] the natural map $\OO_B \to \Eend(\VV^{ab},e^{ab})=\HH^0\EE_\bullet(\VV^{ab},e^{ab})$ is an isomorphism;
	\item[(b)] the Kodaira-Spencer map $KS: \TT_B \to \HH^1\EE_\bullet(\VV^{ab},e^{ab})$ is an isomorphism.
\end{enumerate}
\end{defi}

Condition (a) is satisfied for families of endosimple representations (see \cite[Lem.\ 3.4]{Toda}).
Both conditions are satisfied for the moduli spaces of stable quiver representations corresponding to an indivisible
dimension vector (see \cite[Prop.\ 5.3]{King}).

Let us point out some consequences of the assumptions (a) and (b).
Given $f:S\to B$ (where $S$ is a commutative scheme),
for $(V,e)=(f^*\VV^{ab},f^*e)$ we have 
$$\Eend(V,e)=\HH^0\EE_\bullet(V,e)=\HH^0 Lf^*\EE_\bullet(\VV^{ab},e^{ab})\simeq f^*\HH^0\EE_\bullet(\VV^{ab},e^{ab})\simeq 
f^*\OO_B\simeq\OO_S,$$
where we used the fact that $\HH^1\EE_\bullet(\VV^{ab},e^{ab})\simeq \TT_B$ is locally free.
Also, if $S$ is affine, then
for any coherent sheaf $\FF$ on $S$ we have
$$\HH^1(\EE_\bullet(V,e) \ot \FF)\simeq \HH^1\EE_\bullet(V,e)\ot\FF\simeq f^*\TT_B\ot \FF.$$

Now we consider the following analog of Definition \ref{fun-defi} for quiver representations.

\begin{defi} For an excellent family $(\VV^{ab},e^{ab})$ of representations of $Q$ over a smooth (commutative) base $B$,
we define the functor $h^{NC}_B:\NN\to Sets$ by letting $h^{NC}_B$ to be the set of isomorphism classes of
the following data $(f, V_\La,\phi)$. Let $X=\Spec(\La)$ and let $X^{ab}_0$ be the reduced scheme of the 
abelianization of $X$. Then $f:X^{ab}_0\to B$ is a morphism, $(V_\La,e_\La)$ is a representation of $Q$ over $X$,
and $\phi: (E_{\La},e_\La)|_{X^{ab}_0}\simeq (f^*\VV^{ab},f^*e^{ab}$ is an isomorphism of representations of $Q$.
\end{defi}

We have the following analog of Theorem \ref{loc-rep-thm} (and Proposition \ref{fun-factor-prop}).

\begin{thm}\label{loc-rep-quiver-thm} 
The functor $h^{NC}_B$ is formally smooth and factors through the category $a\NN$.
If the base $B$ is affine then for every $d\ge 0$ the functor $h^{NC}_B|_{a\NN_d}$ is representable by a $d$-smooth
thickening of $B$.
\end{thm}

\Pf . The proof follows the same steps as in the case of families of vector bundles.
We already shown before that $h^{NC}_B$ is formally smooth. The fact that $h^{NC}_B$ factors through $a\NN$ is
proved similarly to Proposition \ref{fun-factor-prop}.

The key technical computation is the analog of Lemma \ref{com-KS-lem}, which in our case claims commutativity of
the diagram
\begin{equation}\label{quiver-KS-comm-sq}
\begin{diagram}
U(f)&\rTo{\De_f}& \Der(A^{ab},I)\\
\dTo{}&&\dTo{-KS}\\
\Aut(\VV_\La,e_\La)&\rTo{\de_0}& H^0(X^{ab},\HH^1\EE_\bullet(\VV^{ab},e^{ab})\otimes \II)
\end{diagram}
\end{equation}
associated with a central extension \eqref{central-ext-eq} and a representation
$(\VV_{\La'},e_{\La'})$ of $Q$ over $X'=\Spec(\La')$. Here we assume that $h_B^{NC}|_{a\NN_{d-1}}$ is represented
by $A\in\NN_{d-1}$, and that $\La\in a\NN_{d-1}$ and $(V_\La,e_\La)$ is a $Q$-representation over 
$X=\Spec(\La)$
corresponding to a homomorphism $f:A\to \La$.
Also, $(\VV_{\La'},e_{\La'})$ is a $Q$-representation over $X'$, extending $(\VV_{\La},e_{\La})$.
The right vertical arrow in \eqref{quiver-KS-comm-sq}
is induced by the KS-map \eqref{quiver-KS-map}, and the bottom arrow
is the connecting map defined in Sec.\ \ref{nonab-hypercoh}. More precisely, we use here the identification for any
quiver representation $(\VV,e)$ over $X$ of the automorphism group
$\Aut(\VV,e)$ with the group
$\hH^0(X,\GG(n_\bullet),c)$, where $c\in \hH^1(X,\GG(n_\bullet)\car\EE_1(n_\bullet))$ is the class of $(\VV,e)$.
Also, we use the natural isomorphism 
\begin{equation}\label{hypercoh-KS-eq}
\hH^1(X,\EE_\bullet(\VV^{ab},e^{ab})\ot \II)\rTo{\sim} H^0(X,\HH^1\EE_\bullet(\VV^{ab},e^{ab})\ot \II)
\end{equation}
induced by the projection $\EE_1(\VV^{ab})\to \HH^1\EE_\bullet(\VV^{ab},e^{ab})$. 

We assume that there is an open covering $(U_i)$ of $B$ and trivializations $\varphi_i^{ab}$ of $\VV^{ab}|_{U_i}$
and the compatible trivializations $\psi_i$ of $\VV_\La$ and $V_{\La'}$ over the covering $\wt{U}_i=q^{-1}U_i$.
Let $(g_{ij}, e_i)$ be the Cech $1$-cocycle corresponding to the universal family over $\Spec(A)$, 
so that the corresponding cocycle for $(\VV_\La,e_\La)$ is $(f(g_{ij}), f(e_i))$.

By definition of $\de_0$ (see Sec.\ \ref{nonab-hypercoh}),
starting from an automorphism $\a$ of $\Aut(\VV_\La,e_\La)$ we can lift it over $\wt{U}_i$ to an automorphism $\a'_i$
of $(\VV_{\La'},e_{\La'})$ and then define $\de_0(\a)$ is the class of the Cech $1$-cocycle with
values in $\EE_\bullet(\VV^{ab},e^{ab})\otimes I$, given by
$$a_{0,ij}=(\a'_i)^{-1}\a'_j-\id, \ \ a_{1,i}=(\a'_i)^{-1}e_{i,\La'}-e_{i,\La'}.$$
Calculating as in the proof of Lemma \ref{com-KS-lem}, and recalling that the action of $\GG_0(n_\bullet)$ on $\EE_1(n_\bullet)$
is given by conjugation, we get
$$a_{0,ij}=\psi_i([\wt{u}^{-1},\wt{f(g_{ij})}]-\id)\psi_i^{-1}=\psi_i \De_f(u^{-1})(f(g_{ij}))f(g_{ij})^{-1}\psi_i^{-1},$$ 
$$a_{1,i}=\psi_i(\wt{u}^{-1}e_{i,\La'}\wt{u}-e_{\La'})\psi_i^{-1}=\psi_i \De_f(u^{-1})(e_i) \psi_i^{-1},$$
where we extend the derivation $\De_f:A\to I$ to matrices with entries in $A$.
Now we note that the image of the class of this Cech $1$-cocycle under the isomorphism \eqref{hypercoh-KS-eq}
is simply the global section of $\HH^1\EE_\bullet(\VV^{ab},e^{ab})\ot\II$ given by
$$(a_{1,i}\mod \im(d_e))=KS(\De_f(u^{-1}))=-KS(\De_f(u)).$$



\begin{thebibliography}{9}
    \bibitem{Brion} M.~Brion, {\it Representations of quivers}, in {\it Geometric methods in representation theory. I}, 103--144, Soc. Math. France, Paris, 2012.
    \bibitem{CQ} J.~Cuntz, D.~Quillen, {\it Algebra extensions and nonsingularity}, J. AMS 8 (1995), 251--289.
    \bibitem{Fedosov} B.V.~Fedosov. {\it A simple geometrical construction of deformation quantization}, J. Diff. Geom. 40 (1994). 213--238.
    \bibitem{Gir} J.~Giraud, {\it Cohomologie non ab\'elienne}, Springer-Verlag (1971).
	\bibitem{HL} D.~Huybrechts, M. Lehn, {\it The Geometry of Moduli Spaces of Sheaves}, Second Edition, Cambridge University Press, (2010).
	\bibitem{Kapranov} M.~Kapranov, {\it Noncommutative geometry based on commutator expansions}, J. Reine Angew. Math. 505 (1998), 73--118.
	\bibitem{KS} M.~Kashiwara, P.~Schapira, {\it Deformation quantization modules}, Ast\'erisque 345, Soc. Math. France (2012).
	\bibitem{King} A.~D.~King, {\it Moduli of representations of finite dimensional algebras}, 
	Quart. J. Math. Oxford Ser. (2) 45 (1994), 515--530.
	\bibitem{Kont} M.~Kontsevich. {\it Deformation quantization of algebraic varieties}, Lett. Math. Phys. 56 (2001). 271--294.
	\bibitem{Manin} Y.I.~Manin, {\it Gauge Field Theory and Complex Geometry}, Springer-Verlag (1991).
	\bibitem{P-BN} A.~Polishchuk, {\it $A_\infty$-structures, Brill-Noether loci and the Fourier-Mukai transform}, Compositio Math. 140 (2004), 459--481.
	\bibitem{PT} A.~Polishchuk, J. Tu, {\it DG-resolutions of NC-smooth thickenings and NC-Fourier-Mukai transforms}, Math. Ann. 360 (2014), 79--156.
	\bibitem{Toda} Y.~Toda, {\it Non-commutative thickening of moduli spaces of stable sheaves}, Compositio Math. 153 (2017), 1153--1195.
	\bibitem{Toda2} Y.~Toda, {\it Non-commutative virtual structure sheaves}, arXiv:1511.00318.
\end{thebibliography}
\end{document}